\documentclass[11pt]{article}
\usepackage{amsmath, amsthm,amsfonts,amssymb,amscd}
\usepackage{color}
\usepackage{geometry}
\usepackage{tikz}
\usepackage{subfigure}
\numberwithin{equation}{section}
\usepackage[T1]{fontenc}
\usepackage{cite}
\usepackage{enumerate}

\usepackage{float}
\usepackage{soul}
\usepackage{marginnote}
\def\disp{\displaystyle}

\def\tto{\;{\lower 1pt \hbox{$\rightarrow$}}\kern -10pt
\hbox{\raise 2pt \hbox{$\rightarrow$}}\;}

\def\hat{\widehat}
\def\Tilde{\widetilde}
\def\tilde{\widetilde}
\def\Bar{\overline}
\def\ra{\rangle}
\def\la{\langle}
\def\ve{\varepsilon}
\def\B{\mathbb B}
\def\IN{\mathbb N}

\def\h{\hfill\Box}
\def\R{\mathbb R}

\def\ox{\bar{x}}

\def\oy{\bar{y}}

\def\ov{\bar{v}}

\def\gph{\mbox{\rm gph}\,}

\def\dom{\mbox{\rm dom}\,}

\def\h{\hfill\triangle}
\def\dn{\downarrow}
\def\O{\Omega}
\def\Lm{\Lambda}

\def\st{\stackrel}
\def\oR{\Bar{\R}}

\def\lm{\lambda}
\def\gg{\gamma}
\def\dd{\delta}
\def\al{\alpha}
\def\kk{\kappa}

\newcounter{lk}

\begin{document}

\newtheorem{Theorem}{Theorem}[section]
\newtheorem{Proposition}[Theorem]{Proposition}
\newtheorem{Remark}[Theorem]{Remark}
\newtheorem{Lemma}[Theorem]{Lemma}
\newtheorem{Corollary}[Theorem]{Corollary}
\newtheorem{Definition}[Theorem]{Definition}
\newtheorem{Example}[Theorem]{Example}
\renewcommand{\theequation}{\thesection.\arabic{equation}}
\normalsize
\def\proof{
\normalfont
\medskip
{\noindent\itshape Proof.\hspace*{6pt}\ignorespaces}}
\def\endproof{$\h$ \vspace*{0.1in}}

\title{\small \bf SECOND ORDER OPTIMALITY CONDITIONS FOR STRONG LOCAL MINIMIZERS VIA SUBGRADIENT GRAPHICAL DERIVATIVE}
\date{}
\author{Nguyen Huy  Chieu\footnote{
 Email: nghuychieu@gmail.com. }, \, \ Le Van  Hien\footnote{Department of Natural Science Teachers, Ha Tinh University,  Ha Tinh, Vietnam; email: lehiendhv@gmail.com.},
\, \  Tran T. A. Nghia\footnote{Department of Mathematics and Statistics, Oakland University, Rochester, MI 48309, USA; email: nttran@oakland.edu. Research of this author was supported by the US National Science Foundation under grant DMS-1816386.},
\,\ Ha  Anh Tuan\footnote{Faculty of Basic Science,  Ho Chi Minh city University of Transport, Ho Chi Minh city,  Vietnam; email:  hatuanhuyhoang@gmail.com}}
\maketitle
{\small \begin{abstract}
This paper is devoted to the study of  second order  optimality conditions for strong local minimizers in the frameworks of unconstrained and constrained optimization problems  in finite dimensions via subgradient graphical derivative. We prove that the positive definiteness of  the subgradient graphical derivative of an extended-real-valued  lower semicontinuous  proper function  at a proximal stationary point  is sufficient  for   the quadratic growth condition. It is also a necessary condition for   the latter property  when  the  function  is either  subdifferentially continuous, prox-regular, twice epi-differentiable or variationally convex. By applying our results to the $\mathcal{C}^2$-cone reducible constrained programs, we establish no-gap second order optimality conditions for (strong) local minimizers  under the  metric subregularity constraint qualification. These results extend the classical  second order  optimality conditions  by surpassing the well-known Robinson's constraint qualification. Our approach also highlights the interconnection between the strong metric subregularity of subdifferential and quadratic growth condition in optimization problems.

\end{abstract}}

\medskip

\normalsize
\section{Introduction}
\setcounter{equation}{0}
In this paper, we mainly study second order  optimality conditions for strong local minimizers to nonsmooth  extended-real-valued  functions with applications to conic programming  in finite dimensional spaces. Strong local minimizer  is  an important concept in optimization at which the {\em quadratic growth condition} is satisfied. In the case of  unconstrained  $\mathcal{C}^2$-smooth optimization  problems, this property of a minimizer is fully characterized by the positive definiteness of the Hessian  of the cost function at stationary points. When the problem is not smooth, several different types of second order {\em directional derivatives} were introduced to study strong minimizers \cite{Au84, BCS99, BS00, RW98,St86,Wa94}. These structures are purely primal due to the involvement  of quantities only on  primal spaces.  Under some constraint qualifications and special regularities, these second order structures could be computed and lead to many optimization applications; see, e.g.,  \cite{BS00, RW98}. 

It is natural to raise the question whether there is any non-primal second order structure that could characterize strong minimizers and if it exists, the concern is about its computability in different frameworks under weaker  regularities.   In 2014,  Arag\'on Artacho and Geoffroy~\cite{AG14}  used {\em  subgradient gradient derivative}, which is a second order construction  as the graphical derivative acting on the subgradient mapping, to depict strong minimizers to convex functions. Due to the involvement of both graphical derivative and subgradient, it is known as a primal-dual structure. The idea of using subgradient graphical derivative to investigate strong local minimizers indeed dates back to Eberhard and Wenczel \cite{EW09}.   Several related second order structures with full computation have been also used in optimization for different purposes in \cite{BLN18,DSZ17,CH17,CHN18,GM17, GO16,RW98}. The approach \cite{AG14} is based on the interconnection between the {\it quadratic growth condition} and the {\em strong metric subregularity} of the convex subdifferential of the cost function studied earlier in \cite{AG08,ZT95}. Without convexity, Drusvyatskiy, Mordukhovich and Nghia \cite{DMN14} showed that the latter property of the limiting subdifferential at  a local minimizer  implies the quadratic growth condition.  The converse implication is also true for a big class of nonconvex semi-algebraic functions \cite{DI15}. It is worth noting that the strong metric subregularity of the subdifferential is a remarkable property in the study of  the linear convergence of various first-order methods \cite{BLN18,DL18, WYYZZ18}. Therefore,   the connection between  the quadratic growth condition and  the strong metric subregularity of the subdifferential is also interesting from numerical viewpoint.


The first aim of this paper is to characterize strong local minimizers of unconstrained nonsmooth optimization problems without convexity  via the subgradient graphical derivative structure. We indeed show in Theorem~\ref{thm0} and ~\ref{thm1} that the positive definiteness of the subgradient graphical derivative at {\em proximal} stationary points is sufficient for strong local minimizers.  It also become necessary  for two broad classes of {\it subdifferentially continuous, prox-regular, twice epi-differentiable functions} \cite[Chapter 13]{RW98} and {\it variationally convex functions} \cite{R18,R18b}.
Through our approach, it is revealed that the appearance of strong local minimizers/quadratic growth conditions to the cost function is equivalent to the strong metric subregularity of subgradient mapping at local minimizers in the latter two favorable settings.   Since functions from the just mentioned two classes are not necessarily semi-algebraic and convex,   our results complement the corresponding ones constructed in \cite{AG08,DMN14,DI15}.

  The second aim of this paper is to study strong local minimizers to smooth conic programming by applying our theory to establishing no-gap second order  optimality conditions for (strong) local minimizers. For nonlinear programming, such a no-gap optimality condition is well-known in literature \cite{BT80,I79} via the  classical  second order necessary and sufficient conditions at stationary points under the Mangasarian-Fromovitz constraint qualification or even under the calmness condition (see \cite[Theorem 2.1]{GLY13} and \cite[Theorem 3.70 (i)]{BS00}).  This result was extended  in \cite{BS00} to the class of $\mathcal{C}^2$-cone reducible programming, which include problems of nonlinear programming, semidefinite programming, and second order cone programming under Robinson's constraint qualification (RCQ). Without RCQ, there are limited results about this important characterization. The  approach in \cite{BS00} is mainly based on some primal structures such as the {\em second order subderivatives} and {\em second order tangent cones}; see also \cite{BCS99}. With the  primal-dual structure, we are able to establish similar results under the so-called {\em metric subregularity constraint qualification} (MSCQ) \cite{GO16,GM17}, which is  strictly weaker than RCQ.  Especially, connection between the quadratic growth condition   and   the strong metric subregularity at local minimizers is also obtained.

The rest of the paper is organized as follows.  Section~2   recalls  some  materials  from variational analysis, which are needed for the sequel analysis.    Section~3 is devoted to the study of the quadratic growth condition for extended-real-valued functions via the subgradient graphical derivative. Section 4 presents our results on no-gap second order necessary and sufficient optimality  conditions for conic programs under the metric subregularity constraint qualification. Finally, section~5 consists of   some concluding remarks on the obtained results as well as the perspectives of this research direction.

\section{Preliminaries}
\setcounter{equation}{0}

In this section we recall some basic notions and  facts  from variational analysis that will be used repeatedly  in the sequel; see  \cite{DR14,M06, M18, RW98} for more details. Let $\O$ be a nonempty subset of the Euclidean space $\R^n$ and $\ox$ be a point in $\O$.  
The (Bouligand-Severi) {\it tangent/contingent cone} to the set $\Omega$ at $\bar x\in \Omega$  is  known as
$$T_\Omega(\bar x):=\big\{v\in\R^n | \, \mbox{there exist}\  t_k \downarrow 0, \  v_k\rightarrow v\ \mbox{ with }\  \bar x+t_kv_k\in\Omega \ \mbox{for all} \   k\in \mathbb{N}\big\}.$$
The  polar cone  of the tangent cone  is  the  (Fr\'echet) {\it regular normal cone} to $\Omega$ at $\bar x$ defined by
 \begin{equation}\label{eqdualTN}\widehat{N}_\Omega(\bar x):=T_\Omega(\bar x)^\circ.\end{equation}
Another normal cone construction used in our work is the  (Mordukhovich) {\it limiting/basic  normal cone} to $\Omega$ at $\bar x\in \Omega$ defined by
$$N_\Omega(\bar x):=\big\{v\in \R^n\, | \,  \mbox{there exist}\,   x_k\st{\Omega}\rightarrow \bar x,\, v_k\in \widehat{N}_\Omega(x_k) \  \mbox{with} \  v_k\rightarrow v\big\}.$$
When $\bar x\not\in \Omega,$ we set $T_\Omega(\bar x)=\emptyset$ and $N_\Omega(\bar x)=\widehat{N}_\Omega(\bar x)=\emptyset$ by convention.
When the set $\Omega$ is convex, the above tangent cone and normal cones  reduce to the tangent cone and normal cone in the  sense of classical convex analysis.

Consider the set-valued mapping $F: \R^n\rightrightarrows\R^m$ with the domain $\dom F:=\big\{x\in \R^n|\; F(x)\neq \emptyset\big\}$ and graph ${\rm gph}\, F:=\big\{(x,y)\in \R^n\times \R^m  |\,   y\in F(x)\big\}$. Suppose that  $(\ox,\oy)$ is an element of ${\rm gph}\, F$.  The {\it graphical derivative} of $F$ at $\bar x$ for $\bar y\in F(\bar x)$ is the set-valued mapping $DF(\bar x|\bar y): \R^n\rightrightarrows\R^m$ defined by
\begin{equation}\label{GraDer}
DF(\bar x|\bar y)(w):=\big\{ z\in\R^m \,  |\,  (w,z)\in T_{{\rm gph}F}(\bar x,\bar y)\big\} \ \,\, \mbox{for }\   w\in\R^n,
\end{equation}
which means   ${\rm gph}\, DF(\bar x|\bar y)= T_{{\rm gph}F}(\bar x,\bar y)$; see, e.g.,  \cite{DR14, RW98}.
   We note further that if $\Phi: \R^n\rightarrow \R^m$ is a single-valued mapping differentiable at $\bar x,$  then $ D\Phi(\bar x|\Phi(\ox))(w)=\nabla \Phi(\bar x)w$ for any $w\in \R^n$.

Following \cite[Section 3.8]{DR14}, we say $F$ is   {\em metrically  subregular} at $\ox\in {\rm dom}\, F$ for $\oy\in F(\ox)$ with modulus $\kk> 0$ if there exists a neighborhood $U$ of $\ox$ such that
\begin{equation}\label{subreg}
d(x;F^{-1}(\oy))\le\kk\, d(\oy;F(x))\quad \mbox{for all}\quad  x\in U,
\end{equation}
where $d(x;\O)$ represents the distance from a point $x\in \R^n$ to a set $\O\subset \R^n$. The infimum of all such $\kk$ is the modulus of metric subregularity, denoted by ${\rm subreg}\, F(\ox|\oy)$. If additionally $\ox$ is an isolated point to $F^{-1}(\oy)$, we say $F$ is   {\em strongly metrically  subregular} at $\ox$ for $\oy$. It is known from \cite[Theorem 4E.1]{DR14} that $F$ is strongly metrically subregular   at $\ox$ for $\oy$ if and only if
\begin{equation}\label{Don}
DF(\ox|\oy)^{-1}(0)=\{0\}.
\end{equation}
Moreover,  in the latter case, its modulus of (strong) metric subregularity  is computed by
\begin{equation}\label{mod}
{\rm subreg}\, F(\ox|\oy)=\dfrac{1}{\inf\{\|z\||\; z\in DF(\ox|\oy)(w), \|w\|=1\}}.
\end{equation}


Assume that  $f: \R^n\rightarrow \overline{\R}:=\R\cup\{\infty\}$ is an extended-real-valued lower semicontinuous (l.s.c.)  proper  function with $\bar x\in\dom f:=\big\{x\in \R^n|\; f(x)<\infty\big\}$. The {\it limiting  subdifferential} (known also as the Mordukhovich/basic subdifferential)  of $f$ at $\bar x$ is defined by
$$\partial f(\bar x):=\big\{ v\in\R^n\, |\, (v,-1)\in N_{{\rm epi} f}(\bar x,f(\ox))\big\},$$
where ${\rm epi}\, f:=\big\{(x,r)\in \R^n\times \R|\;r\ge f(x)\big\}$ is the epigraph of $f$. Another subdifferential construction used in this paper is the {\em proximal subdifferential} of $f$ at $\bar x$ defined by
\begin{equation}\label{Prox}
\partial_p f(\ox):=\Big\{v\in \R^n|\; \liminf_{x\to \ox}\frac{f(x)-f(\ox)-\la v, x-\ox\ra}{\|x-\ox\|^2}>-\infty \Big\}.
\end{equation}
It is well-known that
\begin{eqnarray}\label{LProx}
\partial f(\ox)=\big\{v\in \R^n|\; \mbox{there exists}\; (x_k,v_k)\to (\ox,v)\; \mbox{with}\,v_k\in \partial_p f(x_k)\; \mbox{and}\; f(x_k)\to f(\ox)\big\},
\end{eqnarray}
which shows that $\partial_p f(\ox)\subset \partial f(\ox)$.

  Function $f$ is  said to be {\it prox-regular} at $\ox\in \dom f$ for $\ov\in \partial f(\ox)$ if there exist $r,\ve>0$ such that for all $x,u\in \B_\ve(\ox)$ with $|f(u)-f(\ox)|<\ve$ we have
\begin{equation}\label{prox}
f(x)\ge f(u)+\la v, x-u\ra-\frac{r}{2}\|x-u\|^2\quad \mbox{for all}\quad v\in \partial f(u)\cap \B_\ve(\ov),
\end{equation}
where  $\B_\ve(\ox):=\{x|\;\|x-\ox\|\le\ve\}$ is the closed ball with center $\ox$ and radius $\ve$;    see \cite[Definition~13.27]{RW98}. This  clearly implies that $ \partial f(u)\cap \B_\ve(\ov)\subset \partial_p f(x)$ whenever $\|u-\bar x\|<\ve$ with $|f(u)-f(\ox)|<\ve$.  Moreover, $f$ is said to be {\em subdifferentially continuous} at $\ox$ for $\ov$  if whenever $(x_k, v_k)\to (\ox,\ov)$ and $v_k\in \partial f(x_k)$, one has $f(x_k)\to f(\ox)$; \cite[Definition~13.28]{RW98}. In the case $f$ is  subdifferentially continuous at $\ox$ for $\ov$, the inequality   ``$|f(u)-f(\ox)|<\ve$'' in the definition of prox-regularity above could be omitted.


Recall \cite[Definition 13.3]{RW98} that the {\it second subderivative} of $f$  at $\bar x$ for $v\in \R^n$ and $w\in\R^n$  is given by
\begin{equation}\label{ssd}
d^2f(\bar x|v)(w)=\liminf\limits_{ \begin{subarray}\quad \,\ \tau\searrow 0\\
w'\longrightarrow w\end{subarray}}\Delta^2_\tau f(\bar x|v)(w'),
\end{equation}
where
$$\Delta^2_\tau f(\bar x|v)(w')=\frac{f(\bar x+\tau w')-f(\bar x)-\tau \langle v, w'\rangle}{ \frac{1}{2}\tau^2}.$$
Function  $f$ is said to be {\it twice epi-differentiable} at $\bar x\in \R^n$ for $v\in \R^n$ if
for every $w\in \R^n$ and choice of $\tau_k\searrow 0$ there exist $w^k\to w$ such that
$$\frac{f(\bar x+\tau_k w^k)-f(\bar x)-\tau_k \langle v, w^k\rangle}{ \frac{1}{2}\tau_k^2}\to d^2f(\bar x|v)(w);$$
see, e.g., \cite[Definition 13.6]{RW98}. We note that {\em fully amenable functions},  including the maximum of finitely many $C^2$-functions,  are important examples for    subdifferentially continuous prox-regular and twice epi-differentiable l.s.c. proper functions  \cite[Corollary~13.15 \& Proposition~13.32]{RW98}.

The main second order  structure used in this paper is the {\it subgradient graphical derivative} $D(\partial f)(\ox|\ov): \R^n\rightrightarrows \R^n$ at $\ox$ for $\ov\in \partial f(\ox)$, which  is defined from \eqref{GraDer} by
 \begin{equation}\label{GD}
 D(\partial f)(\ox|\ov)(w):=\big\{z\, |\, (w,z)\in T_{{\rm gph}\,\partial f}(\ox,\ov)\big\}\quad \mbox{for all}\quad w\in \R^n.
 \end{equation}
In the case that $f$ is twice epi-differentiable, prox-regular, subdifferentially continuous at $\ox$ for $\ov$, it is known from  \cite[Theorem~13.40]{RW98} that
\begin{equation}\label{Dh}
D(\partial f)(\ox|\ov)=\partial h\quad \mbox{with}\quad h=\frac{1}{2}d^2f(\bar x|\ov),
\end{equation}
which is an important formula in our study. When $f$ is twice differentiable at $\ox$, it is clear that $D(\partial f)(\ox|\nabla f(\ox))=\nabla^2f(\ox)$.
 \section{Second Order Optimality Conditions for Strong Local Minimizer  via   Subgradient Graphical Derivative}\label{Sec3}
\setcounter{equation}{0}


%
%
Given a function   $f\colon\R^n\to\oR$ and a point  $\ox\in \dom f$,  $\bar x$ is called  a {\it strong local minimizer} of $f$  with modulus $\kk>0$ if there is a number $\gg>0$ such that the following {\it quadratic growth condition}  (QGC, in brief) holds
\begin{equation}\label{GC}
f(x)-f(\ox)\ge\frac{\kk}{2}\|x-\ox\|^2\quad \mbox{for all}\quad x\in \B_\gg(\ox).
\end{equation}
We define the exact modulus for QGC of  $f$ at $\ox$ by
\[
{\rm QG}\, (f;\ox):=\sup\big\{\kk>0|\; \ox \mbox{ is a strong local minimizer of $f$ with modulus $\kk$}\big\}.
\]

In this section we introduce several new  sufficient and necessary   conditions for the quadratic growth condition  \eqref{GC} by using the second order construction defined in \eqref{GD}. The following result taken from \cite[Corollary 3.5]{DMN14} providing a  sufficient condition for the QGC  of $f$ at $\ox$ \eqref{GC} via strong metric subregularity on the subgradient mapping is a significant tool in  our analysis.

\begin{Lemma}{\bf (strong metric subregularity of the subdifferential, \cite[Corollary~3.5]{DMN14}).}\label{DMN14} Let $f\colon \R^n\to\oR$ be a l.s.c. proper function and   let  $\ox\in\dom f$ be a stationary point of $f$   with  $0\in \partial f(\ox)$. Suppose that  the subgradient mapping $\partial f$ is strongly metrically subregular at $\ox$ for $0$ with modulus $\kk>0$ and there are real numbers $r\in(0,\kk^{-1})$ and $\delta>0$ such that
\begin{eqnarray}\label{3.8}
f(x)\ge f(\ox)-\frac{r}{2}\|x-\ox\|^2\quad\mbox{for all}\quad x\in\B_\delta(\ox).
\end{eqnarray}
Then for any $\al\in (0,\kk^{-1})$, there exists a real number $\eta>0$ such that
\begin{eqnarray}\label{3.9}
f(x)\ge f(\ox)+\frac{\al}{2}\|x-\ox\|^2\quad\mbox{for all}\quad x\in\B_\eta(\ox).
\end{eqnarray}
\end{Lemma}

When the function $f$ is convex, the QGC could be fully characterized via the positive definiteness of subgradient graphical derivative \eqref{GD}  \cite[Corollary~3.7]{AG14}. Without convexity, we show in the following result that  such a property is sufficient for  QGC.

\begin{Theorem}{\bf (Sufficient condition  I for strong local minimizers via subgradient graphical derivative).}\label{thm0} Let $f:\R^n\to\oR$ be a proper l.s.c.  function with  $\ox\in \dom f$. Suppose that   $0\in\partial_p f(\ox)$ and that there exists some  real number $c>0$ such that
\begin{equation}\label{SO}
\la z,w\ra\ge c\|w\|^2\quad \mbox{for all}\quad z\in D(\partial f)(\ox|0)(w)\quad \mbox{and}\quad w\in \dom D(\partial f)(\ox|0).
\end{equation}
Then $\ox$ is a strong local minimizer with any modulus $\kk\in (0, c)$. Moreover, we have
\begin{equation}\label{eq}
{\rm QG}(f;\ox)\ge\inf\Big\{\frac{\la z,w\ra}{\|w\|^2}\Big|\; z\in D (\partial f)(\ox|0)(w), w\in \dom D(\partial f)(\ox|0)\Big\}
\end{equation}
with the convention that $0/0=\infty$.
\end{Theorem}
\noindent{\bf Proof.} Since $0\in \partial_p f(\ox)$, we have $0\in \partial f(\ox)$ and there exist $r, \gg>0$ such that
\begin{equation}\label{QB}
f(x)-f(\ox)\ge -\frac{r}{2}\|x-\ox\|^2\quad \mbox{for all}\quad x\in \B_\gg(\ox).
\end{equation}
 To proceed, pick any $s>r$ and define $g(x):= f(x)+\frac{s}{2}\|x-\ox\|^2$, it is clear that
\begin{equation}\label{GC2}
g(x)-g(\ox)\ge  \frac{s}{2}\|x-\ox\|^2\quad \mbox{for all}\quad x\in \B_\gg(\ox).
\end{equation}
Note further that $\partial g(x)=\partial f(x)+s(x-\ox)$ and thus  $0\in\partial g(\ox)$. Thanks to the sum rule of graphical derivative  \cite[Proposition~4A.2]{DR14}, we have
\begin{equation}\label{sume}
D(\partial  g)(\ox|0)(w)=D(\partial f)(\ox|0)(w)+sw\quad \mbox{for all}\quad w\in \R^n.
\end{equation}
Take any $(z,w)\in \R^n\times \R^n$ with $z\in D(\partial  g)(\ox|0)(w)$, i.e., $z-sw\in D(\partial f)(\ox|0)(w)$. It follows from \eqref{SO} that $\la z-sw,w\ra\ge c\|w\|^2$, which means
\begin{equation}\label{mod2}
\|z\|\cdot \|w\|\ge \la z,w\ra\ge (c+s)\|w\|^2.
\end{equation}
We obtain that $D(\partial g)(\ox,0)^{-1}(0)=\{0\}$, i.e., $\partial g$ is strongly metrically  subregular at $\ox$ for $0$ by \eqref{Don}. Moreover, by \eqref{mod} tells us that
\[
{\rm subreg}\, \partial g(\ox|0)\le (c+s)^{-1}.
\]
Since $\ox$ is a local minimizer of $g$ by \eqref{GC2}, it follows from Lemma~\ref{DMN14} that for any $\ve>0$ there exists $\eta\in (0,\gg)$ such that
\begin{equation*}
\begin{array}{rl}
g(x)&\ge g(\ox)+\frac{1}{2({\rm subreg}\, \partial g(\ox|0)+\ve)}\|x-\ox\|^2\\
&\ge  g(\ox)+\frac{1}{2((c+s)^{-1}+\ve)}\|x-\ox\|^2,
\end{array}
\end{equation*}
for all $x\in \B_\eta(\ox).$ Since $f(x)=g(x)-\frac{s}{2}\|x-\ox\|^2$, we obtain from the latter  that
\begin{equation}\label{GC3}
f(x)\ge f(\ox)+\frac{1}{2}\Big[\frac{1}{(c+s)^{-1}+\ve}-s\Big]\|x-\ox\|^2= f(\ox)+\frac{1}{2}\frac{\frac{c}{c+s}-s\ve}{(c+s)^{-1}+\ve}\|x-\ox\|^2.
\end{equation}
By choosing $\ve>0$ sufficiently small, $\ox$ is a strong local minimizer of $f$ with a positive modulus being smaller than but arbitrarily close to $c$. This  verifies that $\ox$ is a strong local minimizer of $f$ with any modulus in $(0, c)$  and the QGC of $f$ holds at $\ox$. Moreover, the inequality \eqref{eq} follows from \eqref{GC3} when taking $\ve\dn 0$ and $c\to$   the infimum on  the right-hand side of \eqref{eq}. The proof is complete. \endproof

When the function $f$ is twice differentiable, the above theorem recovers the classical second order  sufficient condition, which says if $\nabla f(\ox)=0$ and there exists some $c>0$ such that
\begin{equation}\label{Cla}
\la\nabla^2 f(\ox)w,w\ra\ge c\|w\|^2\quad \mbox{for all}\quad w\in \R^n,
\end{equation}
then $\ox$ is a strong local minimizer of $f$. Condition \eqref{Cla} is  known to be  equivalent to the  condition:
\begin{equation}\label{Cla2}
\la\nabla^2 f(\ox)w,w\ra>0\quad \mbox{for all}\quad w\in \R^n, w\neq 0.
\end{equation}
In the nondifferential case, it is natural to question whether condition \eqref{SO} is  equivalent to the following  condition:
\begin{equation}\label{SO2}
\la z,w\ra>0 \quad \mbox{for all}\quad z\in D(\partial f)(\ox|0)(w), w\in \dom D(\partial f)(\ox|0)\setminus \{0\}.
\end{equation}
Obviously, this inequality is a consequence of \eqref{SO}. In the general case, we do not know yet whether the converse  implication is also true. However, we show that \eqref{SO2} is also a sufficient condition to the QGC in the next result, which could be seen as a refinement for our Theorem~\ref{thm0} above. The equivalence between \eqref{SO} and \eqref{SO2} will be clarified later in Theorem~\ref{thm3}, Theorem~\ref{nevc}, and Theorem~\ref{main4} for several broad classes of nondifferentiable functions.

It is worth noting that both \eqref{SO} and \eqref{SO2} imply the strong metric regularity of the subgradient mapping $\partial f$ at $\ox$ for $0$ due to \eqref{Don}. This important feature allows us to use Lemma~\ref{DMN14} to verify strong local minimizer in the following theorem.

\begin{Theorem}{\bf (Sufficient condition II for strong local minimizers via subgradient graphical derivative).}\label{thm1} Let $f:\R^n\to\oR$ be a proper l.s.c.  function with   $\ox\in \dom f$. Consider the following assertions:
\begin{enumerate}[{\bf(i)}]
\item $\ox$ is  a strong local minimizer of $f$.
\item $\ox$ is a local minimizer and $\partial f$ is strongly metrically subregular at $\ox$ for $0$.

\item  $0\in\partial_p f(\ox)$ and  $D(\partial f)(\ox|0)$ is positive definite in the sense of \eqref{SO2}.
\end{enumerate}
Then we have the implications $[{\bf (iii)}\Rightarrow {\bf (ii)}\Rightarrow {\bf (i)}]$.
\end{Theorem}
\noindent{\bf Proof.} The implication $[{\bf (ii)}\Rightarrow {\bf (i)}]$ follows directly from Lemma~\ref{DMN14}. To justify $[{\bf (iii)}\Rightarrow {\bf (ii)}]$, suppose that $0\in \partial_p f(\ox)$ and  condition  \eqref{SO2} is satisfied. It follows  that $$D(\partial f)(\ox|0)^{-1}(0)=\{0\}.$$ By \eqref{Don},   $\partial f$ is strongly metrically subregular at $\ox$ for $0$ with some $\kk>0$. Since $0\in \partial_p f(\ox)$, we find some $r,\gg>0$ such that \eqref{QB} is valid.  Pick any $s>r$ and define $g(x):= f(x)+\frac{s}{2}\|x-\ox\|^2$ as in the proof of Theorem~\ref{thm0} again. For any $(z,w)\in \R^n\times \R^n$ with $z\in D(\partial  g)(\ox|0)(w)$, we derive from \eqref{sume} that  $z-sw\in D(\partial f)(\ox|0)(w)$. It follows  from \eqref{SO2} that $\la z-sw,w\ra\ge 0$, which means
\begin{equation*}
\|z\|\cdot \|w\|\ge \la z,w\ra\ge s\|w\|^2.
\end{equation*}
This together with \eqref{Don} and \eqref{mod} tells us that $\partial g$ is strongly metrically  subregular at $\ox$ for  $0$ with
\[
{\rm subreg}\, \partial g(\ox|0)\le s^{-1}.
\]
Since $\ox$ is a local minimizer of $g$ by \eqref{GC2}, it follows from Lemma~\ref{DMN14} again that for any $\ve>0$ with $\frac{s\ve}{s^{-1}+\ve}<\kk^{-1}$, there exists $\dd\in (0,\gg)$ such that
\begin{equation*}
\begin{array}{rl}
g(x)&\ge g(\ox)+\frac{1}{2({\rm subreg}\, \partial g(\ox|0)+\ve)}\|x-\ox\|^2\\
&\ge  g(\ox)+\frac{1}{2(s^{-1}+\ve)}\|x-\ox\|^2,\end{array}
\end{equation*}
for all $x\in \B_\dd(\ox).$
Since $f(x)=g(x)-\frac{s}{2}\|x-\ox\|^2$, we derive
\begin{equation*}
f(x)\ge f(\ox)+\frac{1}{2}\Big[\frac{1}{s^{-1}+\ve}-s\Big]\|x-\ox\|^2= f(\ox)-\frac{1}{2}\frac{s\ve}{s^{-1}+\ve}\|x-\ox\|^2 \quad \mbox{for all}\quad x\in \B_\dd(\ox).
\end{equation*}
Since $\frac{s\ve}{s^{-1}+\ve}<\kk^{-1}$ and $\partial f$ is strongly metrically subregular at $\ox$ for $0$ with modulus $\kk$, $\ox$ is a (strong) local minimizer of $f$ by Lemma~\ref{DMN14}. The proof is complete. \endproof


As far as we know,   the first  idea of using the subgradient graphical derivative to study the quadratic growth condition was initiated by  Eberhard and  Wenczel \cite{EW09} in which they introduced the so-called {\em sufficient condition of the second kind}.
\begin{Definition}{\bf(sufficient condition of the second kind,\cite{EW09})}\label{ew} {\rm  Let $f:\R^n\to\oR$ be a proper l.s.c.  function with   $\ox\in \dom f$ and  $0\in \partial_p f(\ox)$. We say {\it the sufficient condition of the second kind}  holds at $\ox$ when   there exists $\kk>0$ such that for each $w\in \dom D(\partial_p f)(\ox|0)$ with $\|w\|=1$
\begin{equation}\label{exi}
\exists\, z\in D(\partial_p f)(\ox|0)(w) \; \;{\rm satisfying }\;\; \la z,w\ra\ge \kk.
\end{equation}
}\end{Definition}
Precisely, \cite[Theorem~71(2)]{EW09}  claims  that when the function $f:\R^n \to \oR$ is l.s.c., {\em prox-bounded}, and {\em proximally stable}, the sufficient condition of the second kind at $\ox$ with $0\in \partial_p f(\ox)$ ensures the QGC of $f$ at $\bar x$. However, it seems to us that this result is incorrect even in the convex case. To see this, let us consider the following example.


\begin{Example}\label{mainex} {\rm  Define the function $f:\R\to \R$ by
\begin{eqnarray}
f(x)=\left\{\begin{array}{ll}x&\mbox{if}\quad x>1\\
\al_{n+1}x+\beta_{n+1}\quad &\mbox{if}\quad \al_{n+1}< x\le\al_n, n= 0,1, 2, \ldots\\
\beta\quad &\mbox{if}\quad x=0\\
f(-x) \quad &\mbox{if}\quad x<0,
\end{array}\right.
\end{eqnarray}
where $\al_n=1/(n+1)!$, $\beta_{n+1}:=\disp\sum_{k=0}^{n}\frac{1}{k!(k+2)!}$ with $n=0,1,2, \ldots$, $\beta_0=0$, and $\beta=\disp\lim_{n\to \infty}\beta_n$. It is easy to see that $f$ is  a continuous and  convex function with global optimal solution $\ox:=0$, which clearly implies that $f$ is {\em prox-bounded} and {\em proximally stable} at $\ox$ in the sense of \cite{EW09}.

\begin{figure}[H]
\centering
\includegraphics[scale=0.5]{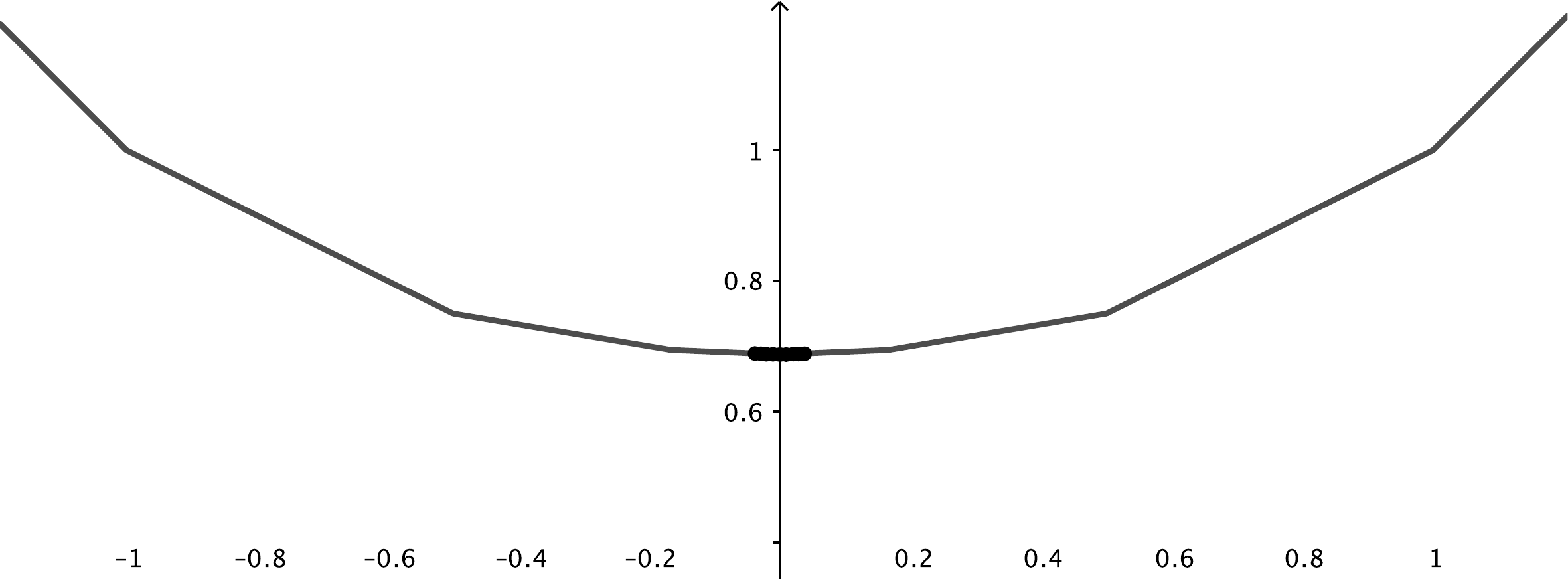}
\caption{Graph of the continuous convex function $f$.}
\end{figure}
\noindent Moreover, direct computation on $\partial_p f$ gives us that
\begin{eqnarray}
\partial_p f(x)=\partial f(x)=\left\{\begin{array}{ll}\{1\}\quad &\mbox{if}\quad x>1\\
 {[}\al_{n+1},\al_n]\quad &\mbox{if}\quad x=\al_n, n=0, 1,2, \ldots\\
\{\al_{n+1}\} \quad &\mbox{if}\quad  \al_{n+1}< x<\al_n, n=0, 1, 2, \ldots\\
\{0\}\quad &\mbox{if}\quad x=0\\
-\partial f(-x) \quad &\mbox{if}\quad x<0.
\end{array}\right.
\end{eqnarray}
\begin{figure}[H]
\centering
\includegraphics[scale=0.35]{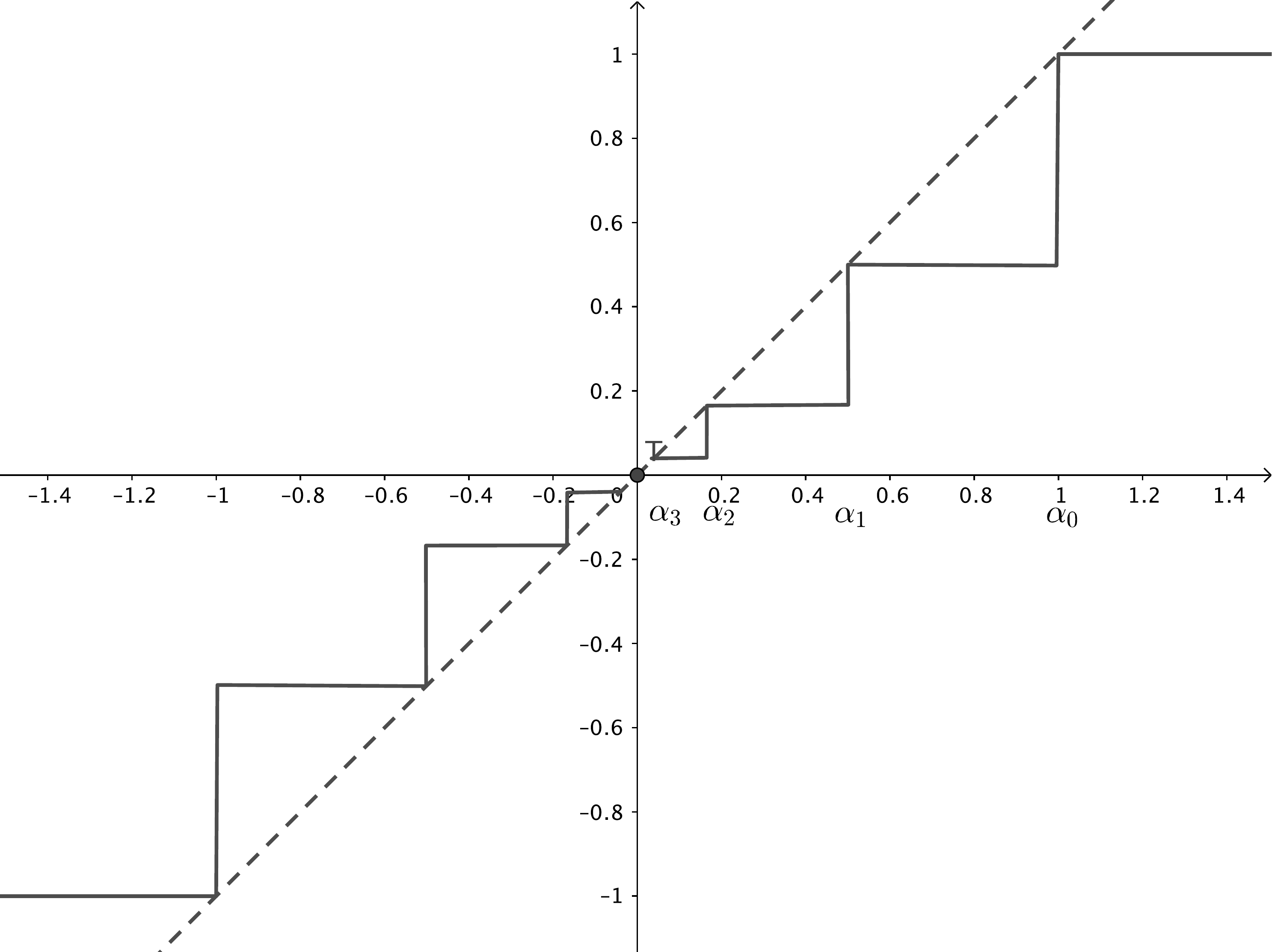}
\caption{Graph of the subgradient mapping $\partial_pf$}
\end{figure}
\noindent Define further $K:=\{(w,z)|\; 0\le z\le w\}\cup\{(w,z)|\; 0\ge z\ge w\}$, we have $\gph \partial_p f\subset K$ and
\begin{equation}\label{tang}
T_{{\rm gph}\, \partial_p f}(\ox,0)\subset T_K(\ox,0)=K.
\end{equation}
 Next we verify the ``$\supset$'' inclusion in \eqref{tang}. Take any $(w,z)$ with $0\le z\le w$ and consider the following three cases:
\begin{itemize}
\item Case 1: {$(w,z)=(0,0)$}  clearly belongs to $T_{{\rm gph}\, \partial_p f}(\ox,0)$.

\item Case 2: $z=0< w$. Choose $t_n= \al_n/w\dn 0$ as $n\to \infty$, we have $t_n(w,\frac{w}{n+2})=(\al_n,\al_{n+1})\in \gph \partial_p f$ and thus $(w,\frac{w}{n+2})\to(w,0)\in T_{{\rm gph}\, \partial_p f}(\ox,0)$.
\item Case 3: $0<z\le w$. Fix $k\in\IN$ satisfying $1/(k+2)\le z/w$ and define $t_n:=\al_n/w$ for $n\ge k$, we have $t_n(w,z)=(\al_n, \al_nz/w)\in\{\al_n\}\times [\al_{n+1},\al_n]  \subset \gph  \partial_p f$.
\end{itemize}
It follows that $\{(w,z)|\; 0\le z\le w\}\subset T_{\gph \partial_p f}(\ox,0)$. Similarly, we have  $\{(w,z)|\; 0\ge z\ge w\}\subset T_{\gph \partial_p f}(\ox,0)$. This together with \eqref{tang} ensures the equality in \eqref{tang}. Note further that   $-1\in D(\partial_p f)(\ox|0)(-1)$ and $1\in D(\partial_p f)(\ox|0)(1)$. Thus the sufficient condition of the second kind \eqref{exi} holds at $\ox$ with $\kk=1$. However, both \eqref{SO} and \eqref{SO2} are not satisfied and  the quadratic growth condition \eqref{GC} is not valid at $\ox$. This tells us that \cite[Theorem~71(2)]{EW09}  is  inaccurate even in the convex case.
}\endproof
\end{Example}
\vspace{-0.1in}
As discussed in the Introduction, QGC and strong local minimizer could be fully characterized via several different types of second order  directional derivatives \cite{Au84, BS00, St86, Wa94,RW98}. For instance, it follows from \cite[Theorem~13.24]{RW98} that $\ox$ is a strong local minimizer to a proper function $f:\R^n \to \R$ if and only if $0\in \partial f(\ox)$ (or $0\in \partial_p f(\ox)$) and the second subderivative \eqref{ssd}  of $f$ at $\ox$ for $0$ is {\em positive definite} in the sense that
\begin{equation}\label{SOD}
d^2f(\ox|0)(w)>0\quad \mbox{for all}\quad w\in \R^n, w\neq 0.
\end{equation}
It is clear that second subderivative \eqref{ssd} is a construction on primal space, while the subgradient graphical derivative \eqref{GD} includes both primal and dual spaces.  Connection between these two constructions could be found in \eqref{Dh} for a special class of   subdifferentially continuous, prox-regular and   twice epi-differentiable functions.  Despite the simplicity of second subderivative and the full characterization of QGC \eqref{SOD}, computing $d^2f(\ox|0)$ could be challenging under some strong regularity conditions. On the other hand,  subgradient graphical derivative is fully computed in many broad classes of optimization problems \cite{CH17, DSZ17, GM17} under milder assumptions.

Unlike \eqref{SOD}, both of our conditions \eqref{SO} and \eqref{SO2} are not  generally necessary conditions for strong local minimizers, as shown in the following example.
\begin{Example}{\rm Let $f:\R\to\R$ be the function defined as follows
$$f(x)=\begin{cases} x\quad\ \ \quad\quad \quad \mbox{if}\  x\in \{0\}\cup [1,+\infty),\\
 \frac{1}{2^n}\quad\, \quad\quad \quad \mbox{if}\ x\in [\frac{3}{2^{n+2}, \frac{1}{2^n})}, \ n=0,1,2,...,\\
2x-\frac{1}{2^{n+1}}\quad \mbox{if}\ x\in [\frac{1}{2^{n+1}}, \frac{3}{2^{n+2}}),\ n=0,1,2,...,\\
f(-x)\quad\quad\ \, \mbox{if} \ x<0.\end{cases}$$
We see that $f(x)\geq f(0)+|x|^2$ for all $x\in [-1,1]$, which means that  $\bar x=0$ is a strong local minimizer of $f$. On the other hand, since $$\left(\bigcup\limits_{n=0}^\infty\big(- \frac{1}{2^n}, -\frac{3}{2^{n+2}}\big)\cup\{0\}\cup \bigcup\limits_{n=0}^\infty\big(\frac{3}{2^{n+2}}, \frac{1}{2^n}\big)\right)\times\{0\}\subset \gph\partial f,$$
it follows that $\R\times\{0\}\subset T_{{\rm gph}\,\partial f}(\bar x,0).$ Therefore, for $ w\in \R\backslash\{0\}$ and $z=0,$ we have
$$z\in D(\partial f)(\bar x|0)(w)\quad\mbox{and}\quad \langle z,w\rangle =0.$$
This shows that \eqref{SO} and \eqref{SO2} are not necessary conditions for strong local minimizers.
}\end{Example}

Our next aim is to present several classes of functions at which both \eqref{SO} and \eqref{SO2} are also  necessary conditions for strong local minimizers. To this end, we first need the following lemma.
\begin{Lemma}\label{lem} Let $h:\R^n \to  \oR$ be  a proper function.  Suppose that $h$ is {\em positively homogenenous of degree} $2$ in the sense that  $h(\lm w)=\lm^2h(w)$ for all $\lm>0$ and $w\in \dom h$.  Then for any $w\in \dom h$ and $z\in \partial h(w)$, we have $\la z,w\ra=2h(w)$.
\end{Lemma}
\noindent{\bf Proof.} For any $z\in \partial h(w)$ with $w\in \dom h$, by \eqref{Prox} and \eqref{LProx} we find  sequences $\{w_k\}\subset\dom h$, $z_k\in \partial_p h(w_k)$, and $\ve_k,r_k>0$ such that $w_k\to w$, $h(w_k)\to h(w)$, $z_k\to z$, and that
\begin{equation*}
h(u)-h(w_k)\ge \la z_k,u-w_k\ra-\frac{r_k}{2}\|u-w_k\|^2\quad \mbox{for all}\quad u\in\B_{\ve_k}(w_k).
\end{equation*}
By choosing $u=\lm w_k\in \B_{\ve_k}(w_k)$ with $0<\lm$ and $\|w_k\|\cdot|\lm-1| <\ve_k$ in the above inequality, the positive homogeneneity of degree $2$ of $h$ tells us that
\begin{equation}\label{inla}
(\lm^2 -1)h(w_k)=h(\lm w_k)-h(w_k)\ge  \la z_k,(\lm-1)w_k\ra-\frac{r_k}{2}(\lm-1)^2\|w_k\|^2.
\end{equation}
When $\lm>1$ satisfying $\|w_k\|\cdot(\lm-1) <\ve_k$,  we get from inequality \eqref{inla} that
\[
(\lm+1)h(w_k)\ge \la z_k, w_k\ra-\frac{r_k}{2}(\lm-1)\|w_k\|^2.
\]
Taking $\lm\dn 1$ gives us that $2h(w_k)\ge\la z_k, w_k\ra$. Similarly, when $0<\lm<1$ with $\|w_k\|\cdot(1-\lm) <\ve_k$, we derive from \eqref{inla} that
\[
(\lm+1)h(w_k)\le \la z_k, w_k\ra-\frac{r_k}{2}(\lm-1)\|w_k\|^2.
\]
By letting $\lm \uparrow 1$, the latter implies $2h(w_k)\le  \la z_k, w_k\ra$. Thus we have $\la z_k, w_k\ra=2h(w_k)$, which clearly yields $\la z,w\ra=2h(w)$ when $k\to \infty$ due to the choice of $z_k, w_k$ at the beginning.   \endproof

\begin{Theorem} {\bf (Characterization of strong local minimizers for prox-regular and   twice epi-differentiable functions).}\label{thm3}  Let  $f: \R^n\to\oR$ be a   l.s.c.  proper function with $\ox\in \dom f$. Suppose that $0\in \partial f(\ox)$ and $f$ is  subdifferentially continuous, prox-regular, and twice epi-differentiable  at $\ox$ for $0$.
Then the following assertions are equivalent:
\begin{enumerate}[{\bf(i)}]

\item  $\bar x$ is a strong local minimizer.

\item $\bar x$ is a  local minimizer and $\partial f$ is strongly metrically subregular at $\bar x$ for $0$.


\item $D(\partial f)(\ox|0)$ is positive definite in the sense of \eqref{SO2}.

\item  $D(\partial f)(\ox|0)$ is positive definite in the sense of \eqref{SO}.

\end{enumerate}

\noindent Moreover,  if one of the assertions ${\bf (i)}-{\bf (iv)}$ holds  then
\begin{equation}\label{eq1}
{\rm QG}(f;\ox)=\inf\Big\{\frac{\la z,w\ra}{\|w\|^2}\Big|\; z\in D (\partial f)(\ox|0)(w)\Big\}.
\end{equation}
\end{Theorem}

\noindent{\bf Proof.}
Since $f$ is subdifferentially continuous and prox-regular at $\ox$ for  $0\in \partial f(\ox)$, we have  $0\in\partial_p f(\ox)$.  Thus, implications  $[{\bf (iv)}\Rightarrow {\bf (iii)}\Rightarrow {\bf (ii)}\Rightarrow {\bf (i)}]$ follow from Theorem~\ref{thm1}.
It remains to   verify $[{\bf (i)}\Rightarrow {\bf (iv)}]$ and \eqref{eq1} is valid.   To this end, suppose that $\bar x$ is a strong local minimizer with modulus~$\kk$ as in \eqref{GC}.  We derive from \eqref{GC} and  \eqref{ssd} that
\begin{equation}\label{eqH0} d^2f(\bar x| 0)(w)\geq \kappa \|w\|^2\quad \mbox{for all}\ w\in \R^n.\end{equation}
Since $f$ is  subdifferentially continuous, prox-regular, and twice epi-differentiable  at $\bar x$ for $0\in \partial f(\bar x),$ it follows from \eqref{Dh}  that
\begin{equation}\label{eqH1}
D(\partial f)(\bar x|0)=\partial h\quad \mbox{with}\quad  h(\cdot):=\dfrac{1}{2} d^2 f(\bar x|0)(\cdot)
\end{equation}
Note from \eqref{ssd} and \eqref{eqH0}  that $h$ is  proper and   positively homogenenous of degree $2$. By Lemma~\ref{lem}, for any $z\in D(\partial f)(\ox|0)(w)=\partial h(w),$  we obtain from \eqref{eqH0} and   \eqref{eqH1} that
\begin{equation}\label{eqH2}
\la z,w\ra=2h(w)=d^2f(\bar x| 0)(w)\geq \kappa \|w\|^2,
\end{equation}
which clearly verifies ({\bf iv}) and $$ \kk\leq \inf\Big\{\frac{\la z,w\ra}{\|w\|^2}\Big|\; z\in D (\partial f) (\ox|0)(w)\Big\}.$$
Since  $\kk$ is an arbitrary modulus of the strong local minimizer $\bar x,$  the latter implies  that
$${\rm QG}(f;\ox)\leq \inf\Big\{\frac{\la z,w\ra}{\|w\|^2}\Big|\; z\in D (\partial f) (\ox|0)(w)\Big\}.$$
This along with \eqref{eq} justifies \eqref{eq1} and finishes the proof.
  \endproof

Besides the full characterization of strong local minimizers in terms of \eqref{SO} and \eqref{SO2} for a class of prox-regular and twice epi-differentiable functions, the above theorem also tells us the equivalence between  QGC and  the strong metric subregularity of subdifferential at a local minimizer for $0$. This correlation has been also established for different classes of functions in \cite{AG08, DMN14, DI15}.

The above theorem allows us to recover \cite[Corollary 73]{EW09}.

\begin{Corollary} Let  $f: \R^n\to\oR$ a   l.s.c. proper   function with $\ox\in \dom f$. Suppose that $0\in \partial f(\ox)$ and that  $f$ is be  subdifferentially continuous, prox-regular ,and twice epi-differentiable  at $\bar x$ for $0.$
Then the following assertions are equivalent:

\begin{enumerate}[{\bf (i)}]
\item  $\bar x$ is a strong local minimizer.

\item The sufficient condition of the second kind  in Definition~\ref{ew} holds at $\ox$.
\end{enumerate}
\end{Corollary}
\noindent{\bf Proof.} Since $f$ is  subdifferentially continuous, prox-regular and twice epi-differentiable  at $\bar x$ for $0\in \partial f(\bar x),$ the proof of Theorem~\ref{thm3}, e.g., \eqref{eqH2} tells us that
\[
\la z,w\ra=d^2f(\bar x| 0)(w)\quad \mbox{for all}\quad  z\in D(\partial f)(\ox|0)(w)=D(\partial_p f)(\ox|0)(w).
\]
Hence, (${\bf ii}$) in this corollary is equivalent to  [(${\bf iv}$), Theorem~\ref{thm3}]. The proof is complete via Theorem~\ref{thm3}.\endproof

The following concept of variational convexity is introduced  recently by Rockafellar \cite[Definition 2]{R18}.

\begin{Definition}{{\bf (variational convexity)}. Let $f: \R^n\to \bar\R$ be a l.s.c.  proper function and  $(\bar x,\bar v)\in \gph\partial f.$ One says that $f$ is  {\it variationally convex} at $\bar x$ for $\bar v$ if there exist an open neighborhood $X\times V$ of $(\bar x,\bar v)$ and a convex lsc function $\hat f\leq f$ on $X$ and $\ve>0$ such that
$$[X_\ve\times V]\cap \gph \partial f= [X\times V]\cap \gph \partial \hat f,$$
and $f(x)=\hat f(x)$ for every $x\in \Pi_X \big([X_\ve\times V]\cap \gph \partial  f\big),$ where
$X_\ve:=\{x\in X | f(x)<f(\bar x)+\ve\}$ and $\Pi_X: \R^n\times \R^n\to \R^n$ is the mapping given by  $\Pi_X(x,v)=x$ for $x\in \R^n$ and $v\in \R^n.$
}
\end{Definition}

The class of variationally convex functions includes convex functions. However,  it  may contain non-convex functions \cite{R18b}.  Note further that  variational convexity implies prox-regularity and subdifferential continuity  \cite{R18b}. The following result resembles Theorem~\ref{thm3} for the class of variationally convex function. 

\begin{Theorem} {\bf (Characterization of strong local minimizer for  variationally convex function).}\label{nevc} Let $f: \R^n\to \bar\R$ be a  l.s.c. proper  function with $\ox\in \dom f$. Suppose that  $0\in \partial f(\bar x)$ and  that $f$ is   variationally convex at $\bar x$ for~$0.$  Then the following assertions are equivalent:

\begin{enumerate}[{\bf (i)}]

\item  $\bar x$ is a strong local minimizer.

\item $\bar x$ is a  local minimizer and $\partial f$ is strongly metrically subregular at $\bar x$ for $0$.

\item $D(\partial f)(\ox|0)$ is positive definite in the sense of \eqref{SO2}.

\item  $D(\partial f)(\ox|0)$ is positive definite in the sense of \eqref{SO}.

\end{enumerate}
\noindent Moreover,  if one of the assertions ${\bf (i)}-{\bf (iv)}$ holds  then
\begin{eqnarray}\label{eq2}
{\rm QG}(f;\ox)\ge \inf\Big\{\frac{\la z,w\ra}{\|w\|^2}\Big|\; z\in D (\partial f)(\ox|0)(w)\Big\}\ge \frac{1}{2} {\rm QG}(f;\ox).
\end{eqnarray}

\end{Theorem}
\noindent{\bf Proof.} Note from  the variational convexity of $f$ at $\bar x$ for $0\in \partial f(\bar x)$  and \eqref{prox} that
$0\in \partial \hat f(\bar x)$ and $0\in \partial_p f(\ox)$.  Similarly to the proof of Theorem~\ref{thm3}, we only need to  verify $[{\rm \bf (i)\Rightarrow(iv)}]$ and the right inequality in \eqref{eq2} due to \eqref{eq}. Suppose that $\bar x$ is a strong local minimizer with modulus $\kk,$ that is,  there is a number $\gg>0$ such that
\begin{equation}\label{SLM}
f(x)-f(\ox)\ge\frac{\kk}{2}\|x-\ox\|^2\quad \mbox{for all}\quad x\in \B_\gg(\ox).
\end{equation}
The variational convexity  of  $f$  at $\bar x$ for~$0$ allows us to find an open neighborhood $X\times V$ of $(\bar x,\bar v)$ and a convex lsc function $\hat f\leq f$ on $X$  with
$$[X\times V]\cap \gph \partial f= [X\times V]\cap \gph \partial \hat f,$$
and  $f(x)=\hat f(x)$ for every $x\in \Pi_X \big([X\times V]\cap \gph \partial  f\big).$
Pick any $z\in D\partial f(\ox|0)(w)$, by \eqref{GD} there exist $t_k\downarrow0$ and  $(z_k,w_k)\to (z,w)$ such that
 $$(\bar x, 0)+t_k(w_k,z_k)\in [X\times V]\cap \gph \partial f= [X\times V]\cap \gph \partial \hat f.$$
 Note that $\bar x, \bar x+t_kw_k\in \Pi_X \big([X\times V]\cap \gph \partial  f\big).$ It follows from \eqref{SLM} that
\begin{eqnarray*}\hat f(\bar x+t_kw_k)=f(\bar x+t_kw_k)\geq f(\bar x)+\frac{\kk}{2}t_k^2\|w_k\|^2= \hat f(\bar x)+\frac{\kk}{2}t_k^2\|w_k\|^2\quad \mbox{for all}\  k\  \mbox{sufficiently large.}
\end{eqnarray*}
Furthermore, since $\hat f$ is convex and  $(\bar x, 0)+t_k(w_k,z_k)\in \gph \partial \hat f,$ we have
$$\hat f(\bar x)-\hat f(\bar x+t_kw_k)\geq -\langle t_kz_k, t_kw_k\rangle\quad \mbox{for all}\  k.$$
Combining the above two inequalities gives us that
$\langle z_k, w_k\rangle \geq  \frac{\kk}{2}\|w_k\|^2$  for sufficiently large $k$ . Letting $k\to \infty$ we get
$\langle z, w\rangle \geq  \frac{\kk}{2}\|w\|^2,$
which clearly clarifies $[{\rm \bf (i)\Rightarrow(iv)}]$  and the right inequality in \eqref{eq2}. \endproof

\section{Second Order Optimality
Conditions for $\mathcal{C}^2$-Reducible Conic Programs under Metric Subregularity Constraint Qualification}
In this section, let us consider the following constrained optimization problem:
  \begin{equation}\label{conicP}  \begin{array}{rl}
   (P)\qquad &\min\limits_{x\in \R^n} \quad g(x)\quad \mbox{subject to} \quad q(x)\in \Theta,
 \end{array}
 \end{equation}
 where $g:\R^n\to \R$ and $q:\R^n\to \R^m$ with $q(x)=\big(q_1(x),...,q_m(x)\big)$ are twice continuously differentiable, and $\Theta$ is a nonempty closed convex subset of $\R^m.$

Define $\Gamma:=\{x\in \R^n\ |\ q(x)\in \Theta\}$ as the feasible solution set to problem \eqref{conicP} and  fix $\bar x\in \Gamma$ with $\bar y:=q(\bar x).$ Put
 \begin{equation}\label{funf}
 f(x):=g(x)+\delta_\Gamma(x)\ \mbox{for all}\  x\in \R^n,
 \end{equation}
 where $\delta_\Gamma(x)$ is the indicator function to $\Gamma$, which equals to $0$ when $x\in \Gamma$ and $\infty$ otherwise.
 Problem \eqref{conicP} can be rewritten  as an unconstrained optimization problem:
 $$\min\limits_{x\in \R^n} f(x).$$

 The given point $\bar x\in \Gamma$ is said to be  a {\it strong local minimizer}  to  problem \eqref{conicP} if  there exist  numbers $\kappa>0, \gg>0$ such that
\begin{equation}\label{C2GC}g(x)\geq g(\bar x)+\frac{\kappa}{2}\|x-\bar x\|^2\quad \mbox{for all}\ x\in \Gamma\cap\B_\gg(\bar x), \end{equation}
that is, $\bar x$ is a strong local minimizer of the function $f$ defined above. In the case of \eqref{C2GC}, we say the quadratic growth condition holds at $\ox$ to problem \eqref{conicP}. Moreover, the exact modulus for the problem \eqref{conicP} at $\ox$ denoted by $QG((P);\ox)$ is the supremum of all $\kk$ such that \eqref{C2GC} is satisfied, i.e., $QG((P);\ox)=QG(f;\ox)$.

 We call $\bar x\in \Gamma $ a  {\it local minimizer}  to  \eqref{conicP}  if  condition \eqref{C2GC} holds with some $\gg>0$ and $\kappa=0$. Furthermore, $\ox$ is  a {\em stationary point} when there exists a Lagrange multiplier $\lm\in N_\Theta(q(\ox))$ such that
\begin{equation}\label{sta}
0=\nabla g(\ox)+\nabla q(\ox)^T\lm.
\end{equation}
The set of Lagrange multipliers satisfying \eqref{sta} is denoted by $\Lm(\ox)$.

In what follows, we always assume that  the closed convex set $\Theta$ is $C^2$-{\em cone reducible} at $\bar y$ to a {\em pointed} closed convex cone $C\subset \R^l$ in the sense that there exists a neighbourhood $V\subset \R^m$ of $\bar y$ and a twice continuously differentiable mapping $h: V\to \R^l$ such that
  \begin{equation}\label{hV}
  h(\bar y)=0,\ \nabla h(\bar y) \ \mbox{is surjective, and}\ \Theta\cap V=\big\{y\in V\ |\ h(y)\in C\big\}.
  \end{equation}
It is worth noting that the assumption of reducible sets allows us to cover wide range of optimization problems including nonlinear programming, semidefinite programming, and second order  cone programming; see, e.g., \cite[Section~3.4.4]{BS00}.

  Furthermore, we assume  that the {\em  metric subregularity constraint qualification} (MSCQ, in brief) \cite{GM15,GM17} holds at $\bar x$, which means the set-valued mapping $F(x):= q(x)-\Theta$, $x\in \R^n$ is metrically subregular at $\ox$ for $0$. This condition is well-known to be stable around $\ox$ and strictly weaker than the notable Robinson's constraint qualification (RCQ) at $\ox$:
\begin{equation}\label{RCQ}
0\in {\rm int}\, \{q(\ox)+\nabla q(\ox)\R^n-\Theta\}.
\end{equation}
Furthermore, if MSCQ is satisfied at $\ox$, there exists some $\eta>0$ such that the normal cone $N_\Gamma(x)$ is identical with $\partial_p \delta_\Gamma(x)$ for all $x\in \B_\eta(\ox)$  and is computed by
 \begin{equation}\label{norma}
 N_\Gamma(\ox)=\partial_p\delta_\Gamma(\ox)=\{\nabla q(\ox)^T\lm|\; \lm\in N_\Theta(q(\ox))\}\quad \mbox{for}\quad x\in \B_\eta(\ox).
 \end{equation}

The main purpose of this section is to establish the necessary and sufficient second order  conditions for strong local minimizer under the MSCQ by using the theory developed in Section~\ref{Sec3}. In order to do so, we need to calculate the subgradient graphical derivative on the function $f$ defined in \eqref{funf} at the given point $\ox$ for $0$ due to Theorem~\ref{thm1} under the assumption that $0\in \partial_p f(\ox)$.  Under MSCQ at $\ox$, it follows from \eqref{norma} that
\begin{equation}\label{sumr}
\partial_p f(x)=\nabla g(x)+\partial_p \delta_\Gamma(x)=\partial f(x)\quad \mbox{for}\quad x\in \B_\eta(\ox)
\end{equation}
with the same $\eta$ in \eqref{norma}. By \eqref{sumr},  observe that $0\in \partial_p f(\ox)$ iff $\ox$ is a stationary point to  \eqref{sta}. Furthermore, we have
\begin{equation}\label{grasum}
D(\partial f)(\ox|0)(w)=D(\nabla g+N_\Gamma)(\ox|0)(w)=\nabla^2 g(\ox)w+DN_\Gamma(\ox|-\nabla g(\ox))(w).
\end{equation}
The following result taken from \cite[Corollary 5.4]{GM17} is helpful in our study.

\begin{Lemma}\label{GM17}{\rm(\cite[Corollary 5.4]{GM17}).}  Let $\ox$ be a stationary point to problem \eqref{conicP}. Suppose that MSCQ is satisfied at $\ox$ and that the set $\Theta$ is $\mathcal{C}^2$-reducible at $\oy=q(\ox)$. Then  we have
\begin{equation}\label{DNG}
D N_\Gamma(\bar x,-\nabla g(\bar x))(w)=\left\{\big(\nabla^2(\lambda^Tq)(\bar x)+\widetilde{\mathcal{H}}_\lambda\big)w\ |\ \lambda\in \overline{\Lambda}(\ox;w) \right\}+N_{\overline{K}}(w)\;\; \mbox{for}\;\;w\in \R^n,
\end{equation}
  where $\overline{K}$ is the critical cone defined by
  \begin{equation}\label{critco}
  \overline{K}:=K\big(\bar x, -\nabla g(\bar x)\big):=T_\Gamma(\ox)\cap\{-\nabla g(\ox)\}^\perp
  \end{equation}
  and the set $\overline{\Lambda}(\ox;w)$ is written by
 \begin{equation}\label{Balm}
 \overline{\Lambda}(\ox;w):={\rm argmax}\big\{w^T\nabla^2(\lambda^Tq)(\ox)w+w^T\widetilde{\mathcal{H}}_\lambda w |\  \lambda \in \Lambda\big(\bar x)\big)\big\}\neq\emptyset
 \end{equation}
  with
  \begin{equation}\label{Hti}
  \widetilde{\mathcal{H}}_\lambda w:=\nabla q(\bar x)^T\nabla^2\left\langle \Big(\nabla h\big(q(\bar x)\big)^*\Big)^{-1}(\lambda), h(\cdot)\right\rangle \big(q(\bar x)\big)\nabla q(\bar x)w.
  \end{equation}
\end{Lemma}
Pick any $z\in D(\partial f)(\ox|0)(w)$,  it follows from \eqref{grasum} and \eqref{DNG} that $w\in \dom D(\partial f)(\ox|0)= \overline{K}$ and there exists $\bar \lm\in \overline{\Lm}(\ox;w)$ and $u\in N_{\overline{K}}(w)$ such that
\[
z= \nabla^2 g(\ox)w+\big(\nabla^2(\bar\lambda^Tq)(\bar x)+\widetilde{\mathcal{H}}_{\bar \lambda}\big)w+u.
\]
Hence, we have
\[
\la z,w\ra= \la \nabla^2 g(\ox)w,w\ra+w^T\nabla^2(\bar \lambda^Tq)w+w^T\widetilde{\mathcal{H}}_{\bar \lambda} w+\la u,w\ra.
\]
Since $\overline{K}$ is a cone, $\la u,w\ra=0$. We derive from the latter and \eqref{Balm} that
\begin{eqnarray}\label{wee}
\la z,w\ra=\max \big\{w^T\nabla^2L(\bar x, \lambda)w+w^T\widetilde{\mathcal{H}}_\lambda w\ |\  \lambda \in  \Lambda(\bar x)\big\} \; \mbox{for}\; z\in D(\partial f)(\ox|0)(w), w\in \overline{K}
\end{eqnarray}
with the Lagrange function $L(x,\lm):=g(x)+\la \lm,q(x)\ra$ for $x\in \R^n, \lm\in \R^m$. This together with Theorem \ref{thm1} tells us that the following second order  condition
\begin{equation}\label{seconic}
\max \big\{w^T\nabla^2L(\bar x, \lambda)w+w^T\widetilde{\mathcal{H}}_\lambda w\ |\  \lambda \in  \Lambda(\bar x)\big\}>0\quad \mbox{for all}\quad w\in \overline{K}
\end{equation}
is sufficient for strong local minimizer $\ox$, when  MSCQ holds at the stationary point $\ox$. This is a classical fact \cite[Theorem 3.137]{BS00} established  under the strictly stronger RCQ. Moreover, \eqref{seconic} is  necessary for the strong local minimizer $\ox$ under RCQ. In order to show this fact under MSCQ, we need a few lemmas in preparation.

Let $Q: h^{-1}(V)\to \R^\ell$ be the mapping  given by $Q(x)=h\circ q(x)$ with $h$ and $V$ taken from \eqref{hV}.

\begin{Lemma} {\rm (\cite[Propostion 3.1]{GM17}).} \label{lm-sed02.1} Suppose that MSCQ holds for the system $q(x)\in \Theta$ at $\bar x\in \Gamma$. Then with any $w\in  \R^n$ satisfying $\nabla Q(\bar x)w\in C$ one can find a positive number $\kappa$ such that for any $(z,y)\in \R^n\times \R^{\ell}$ with
  $$\nabla Q(\bar x)z+\langle w, \nabla^2Q(\bar x)w\rangle+y \in T_{C}(\nabla Q(\bar x)w),$$
  there exists $\tilde z\in \R^n$ satisfying the conditions
  $$\nabla Q(\bar x)\tilde z+\langle w, \nabla^2Q(\bar x)w\rangle \in T_{C}(\nabla Q(\bar x)w)\ \mbox{and}\ \|\tilde z-z\|\leq \kappa\|y\|.$$
  The latter condition can be reformulated as the upper Lipschitzian property
 $$\Psi(y)\subset \Psi(0)+\kappa \|y\|\B_{\R^{n}} \quad \mbox{for all}\ y\in \R^{\ell},$$
 where
 $$\Psi(y):=\big\{z\in \R^n \ | \   \nabla Q(\bar x)z+\langle w, \nabla^2Q(\bar x)w\rangle+y \in T_{C}(\nabla Q(\bar x)w)\big\}.$$
  \end{Lemma}

\begin{Lemma}{\rm (see \cite[p.~242]{BS00}).} \label{lm-sed02} If $w\in \R^n$ with $\nabla Q(\bar x)w\in C$  then
 $$ T^2_{C}\big(Q(\bar x), \nabla Q(\bar x) w\big)=T_{C}\big(\nabla Q(\bar x)w\big),$$
 where
$$T^2_{C}\big(Q(\bar x), \nabla Q(\bar x) w\big):=\Big\{ v \, |\,  \exists t_k\searrow 0 \ {\rm s.t.}\ d\big(Q(\ox)+t_kw+\frac{1}{2}t^2_kv;C\big)=o(t_k^2)\Big\}$$
is the outer second order tangent cone to the set $C$  at $Q(\bar x)$ and in the direction $w\in T_C(Q(\ox))$.
  \end{Lemma}

\begin{Lemma}\label{lm-sed01} Let $\bar x$ be a local minimizer of  \eqref{conicP} satisfying MSCQ. Then, for each $w\in K\big(\bar x, -\nabla g(\bar x)\big)$ and $z\in \R^n$  with
\begin{equation}\label{Qz}
\nabla Q(\bar x) z+\langle w, \nabla^2Q(\bar x)w\rangle \in T^2_{C}\big(Q(\bar x), \nabla Q(\bar x) w\big),
\end{equation}
we have
$$\nabla g(\bar x) z + \langle w, \nabla^2g(\bar x)w\rangle\geq 0.$$
\end{Lemma}
\noindent{\bf Proof.}
Take any  $w\in K\big(\bar x, -\nabla g(\bar x)\big)$ and $z\in \R^n$  satisfying \eqref{Qz}. Then there exists $t_k\searrow 0$ such that
\begin{equation}\label{thmSedEq1}
d\Big(Q(\bar x)+t_k\nabla Q(\bar x)w+\frac{1}{2}t^2_k\big(\nabla Q(\bar x) z+\langle w, \nabla^2Q(\bar x)w\rangle\big);C\Big)=o(t_k^2).
\end{equation}
Defining $x(t_k):=\bar x+t_kw+\frac{1}{2}t_k^2z$, we have
$$d\Big(Q\big(x(t_k)\big); C\Big)=d\Big(Q(\bar x)+t_k\nabla Q(\bar x)w+\frac{1}{2}t^2_k\big(\nabla Q(\bar x) z+\langle w, \nabla^2Q(\bar x)w\rangle\big)+o(t_k^2), C \Big)=o(t_k^2).$$
Since MSCQ holds at $\bar x,$  the set-valued mapping $M_Q(x)=Q(x)-C$ is metrically subregular at $\bar x$ for $0=Q(\ox)\in \R^\ell$ by \cite[Lemma 5.2]{GM17}. Then there exist a neighborhood $U$ of $\bar x$ and a real number $\kappa >0$ such that
$$d\big(x;\Gamma\cap h^{-1}(V)\big)= d(x; M_Q^{-1}(0)\big)\leq \kappa d\big(0; M_Q(x)\big)=\kappa d\big(Q(x);C\big)\quad \mbox{for  all}\ x\in U.$$
We assume without loss of generality that $x(t_k)\in U$ for all $k.$ Consequently,
$$d\big(x(t_k);\Gamma\cap h^{-1}(V)\big)\leq \kappa d\Big(Q\big(x(t_k)\big); C\Big)=o(t^2_k).$$
So for each $k$ one can find $\tilde x (t_k)\in \Gamma\cap h^{-1}(V)$ with $\|x(t_k)-\tilde x (t_k)\|=o(t^2_k).$
Noting that $x(t_k)\to \bar x$ as $k\to \infty$ and $\bar x$ is a local optimal solution to \eqref{conicP}, we may assume that
$g\big(\tilde x (t_k)\big)\geq g(\bar x)$ for all $k.$
On the other hand, by the Taylor expansion,
$$ g\big(x(t_k)\big)=g(\bar x)+t_k\nabla g(\bar x)w+\frac{1}{2}t^2_k\big(\nabla g(\bar x) z+\langle w, \nabla^2g(\bar x)w\rangle\big)+o(t_k^2),$$
and  $ g\big(x(t_k)\big)= g\big(\tilde x(t_k)\big)+o(t_k^2).$
We have
$$g(\bar x)+t_k\nabla g(\bar x)w+\frac{1}{2}t^2_k\big(\nabla g(\bar x) z+\langle w, \nabla^2g(\bar x)w\rangle\big)+o(t_k^2)\geq g(\bar x)\ \mbox{for all}\ k.$$
Combining this with $\nabla g(\bar x)w=0$ yields
$$\frac{1}{2}t^2_k\big(\nabla g(\bar x) z+\langle w, \nabla^2g(\bar x)w\rangle\big)+o(t_k^2)\geq 0.$$
This implies that  $\nabla g(\bar x) z+\langle w, \nabla^2g(\bar x)w\rangle\geq 0.$
 \endproof


We are ready to arrive at the first main result of this section providing  second order  necessary optimality conditions for  conic programs with $C^2$-cone reducibility constraint.

   \begin{Theorem}\label{thm1ConicP} {\rm \bf  (second order  necessary optimality conditions for local minimizers  of   conic program  under MSCQ and $\mathcal{C}^2$-cone reducibility constraint).}   Suppose that $\bar x$ is a stationary point of  \eqref{conicP} at which MSCQ is satisfied and that $\Theta$ is $\mathcal{C}^2$-cone reducible at $\oy=h(\ox)$ to a pointed closed convex cone $C$. Consider the following assertions:

   ${\bf (i)}$ $\bar x$   is a local minimizer for   \eqref{conicP}.

 ${\bf (ii)}$ $D(\partial f)(\ox|0)$ is positive semidefinite in the sense that
\begin{equation}\label{SO2C}
\la z,w\ra\geq0 \quad \mbox{for all}\quad z\in D(\partial f)(\ox|0)(w), w\in {\rm dom}\,D(\partial f)(\ox|0).
\end{equation}

${\bf (iii)}$ For each  $w\in {\rm dom}\,D(\partial f)(\ox|0)$ there exists  $z\in D(\partial f)(\ox|0)(w)$ such that
\begin{equation}\label{SO2t}
\la z,w\ra\geq 0.
\end{equation}


${\bf (iv)}$  For each  $w\in \overline{K}$  one has
 $$\max \big\{w^T\nabla^2L(\bar x, \lambda)w+w^T\widetilde{\mathcal{H}}_\lambda w\ |\  \lambda \in  \Lambda\big(\bar x\big)\big\}\geq 0$$
with $\Tilde H_\lm$ defined in \eqref{Hti}.

Then we have  $[{\bf (i)}\Rightarrow {\bf (ii)}\Leftrightarrow {\bf(iii)}\Leftrightarrow {\bf(iv)}].$
    \end{Theorem}
      \noindent{\bf Proof.}  $[{\bf (i)}\Rightarrow {\bf (iv)}]:$ Suppose $\bar x$  is a local minimizer for   \eqref{conicP}.
        For each $w\in \overline{K}$ let us consider the linear conic problem  $(\widetilde P)$ defined as follows
 \[
 (\widetilde P)\quad \begin{array}{rl}& \inf\limits_{z\in \R^n}\quad  \nabla g(\bar x)z+\langle w,\nabla^2 g(\bar x)w\rangle\\
              & \mbox{s.t.}\quad \  \nabla Q(\bar x)z+\langle w, \nabla^2Q(\bar x)w\rangle \in T_C\big(\nabla Q(\bar x)w\big),
      \end{array}
 \]
and its parametric dual   $(\widetilde D)$ given by
\[
(\widetilde D)\quad \begin{array}{rl}& \sup\limits_{\mu\in \R^\ell}  \quad w^T\nabla^2_x\mathcal{L}(\bar x, \mu)w \\
&\mbox{s.t.}\quad\    \mu \in  N_{C}\big(\nabla Q(\bar x)w\big),  \nabla Q(\bar x)^T\mu=-\nabla g(\bar x),\end{array}
\]
 where $\mathcal{L}(x,\mu)=g(x)+\langle \mu, Q(x)\rangle,$ $(x,\mu)\in \R^n\times \R^\ell.$
      By Lemma~\ref{lm-sed02.1},  the feasible set of  $(\widetilde P)$ is nonempty.      Moreover,  by Lemmas~\ref{lm-sed02}\&\ref{lm-sed01}, ${\rm val}(\widetilde P)$  is finite and ${\rm val}(\widetilde P)\geq 0,$ where ${\rm val}(\widetilde P)$  is the optimal value of  $(\widetilde P).$  Due to the upper Lipschitz continuity of the mapping $\Psi$ in Lemma~\ref{lm-sed02.1}, we derive from  \cite[Propositions~2.147 \& 2.186]{BS00}  that
      ${\rm val}(\widetilde P)={\rm val}(\widetilde D)$
      and that the optimal solution set to  $(\widetilde D)$ is nonempty.
With $\lambda=\nabla h\big(q(\bar x)\big)^T\mu$, note   that
 \[\begin{array}{rl} w^T\nabla^2_x\mathcal{L}(\bar x, \mu)w&=w^T\nabla^2_xL(\bar x, \lambda)w +\Big\langle \mu, (\nabla q(\bar x)w)^T\nabla^2h\big(q(\bar x)\big)(\nabla q(\bar x)w)\Big\rangle \\
 &= w^T\nabla^2_xL(\bar x, \lambda)w +\Big\langle \Big(\nabla h\big(q(\bar x)\big)^*\Big)^{-1}(\lambda), (\nabla q(\bar x)w)^T\nabla^2h\big(q(\bar x)\big)(\nabla q(\bar x)w)\Big\rangle \\
 &=w^T\nabla^2_xL(\bar x, \lambda)w+w^T\widetilde{\mathcal{H}}_\lambda w.\end{array}
 \]
We claim further that
\begin{equation}\label{seteq}
\{\lm=\nabla h\big(q(\bar x)\big)^T\mu|\; \mu \in  N_{C}\big(\nabla Q(\bar x)w\big),  \nabla Q(\bar x)^T\mu=-\nabla g(\bar x)\}= \Lm(\ox).
\end{equation}
To justify the ``$\subset$ inclusion, pick any $\mu\in  N_{C}\big(\nabla Q(\bar x)w\big)\subset C^*$ with   $\nabla Q(\bar x)^T\mu=-\nabla g(\bar x)$ and   define $\lm:=\nabla h\big(q(\bar x)\big)^T\mu$. Since $C$ is a convex cone,  we have
\[
\lm= \nabla h\big(q(\bar x)\big)^T \mu\in \nabla h\big(q(\bar x)\big)^T C^*=\nabla h\big(q(\bar x)\big)^TN_C(Q(\ox))=N_\Theta(q(\ox)).
\]
Moreover, it is clear that
\begin{equation}\label{nons}
\nabla q(\ox)^T\lm=\nabla q(\ox)^T\nabla h\big(q(\bar x)\big)^T \mu=\nabla Q(\ox)^T\mu=-\nabla g(\ox).
\end{equation}
This shows that $\lm \in \Lm(\ox)$ and thus verifies the ``$\subset$'' inclusion in \eqref{seteq}.

 To ensure the opposite inclusion in \eqref{seteq}, take any $\lm\in \Lm(\ox)$. Since $\lm\in N_\Theta(q(\ox))=\nabla h(q(\ox))^TN_C(Q(\ox)$, we find $\mu\in N_C(Q(\ox)=C^*$ with $\lm= \nabla h(q(\ox))^T\mu$. It is similar to \eqref{nons} that $\nabla Q(\ox)^T\mu=-\nabla g(\ox)$. Moreover, note from the fact $w\in \overline{K}$ that
 \[
 \la \mu, \nabla Q(\ox)w\ra=\la \nabla Q(\ox)^T\mu, w\ra=\la -\nabla g(\ox),w\ra=0
 \]
 Since $C$ is a convex cone, we get from the latter that$\mu\in N_C(\nabla Q(\ox)w)$. This clearly verifies the ``$\supset$'' inclusion in \eqref{seteq}.

Note further from \eqref{hV} that $\nabla h(q(\ox)^T$ is injective, we derive that  \begin{equation}\label{sup-max1}
            {\rm val}(\widetilde D)=\max \big\{w^T\nabla^2L(\bar x, \lambda)w+w^T\widetilde{\mathcal{H}}_\lm w\ |\  \lambda \in  \Lambda(\bar x)\big\}.
\end{equation}
Since ${\rm val}(\widetilde D)={\rm val}(\widetilde P) \ge 0$,  the latter implies  the assertion {\bf (iv)}.


$[{\bf (iv)}\Leftrightarrow {\bf (iii)}\Leftrightarrow {\bf (ii)}]:$  Take any  $z\in D(\partial f)(\ox|0)(w)$ with $w\in {\rm dom}\,D(\partial f)(\ox|0)$. We obtain  from \eqref{grasum} and  Lemma~\ref{GM17} that ${\rm dom}\,D\partial f(\ox|0)=\overline{K}$. The equivalence $[{\bf (iv)}\Leftrightarrow {\bf (iii)}\Leftrightarrow {\bf (ii)}]$ simply follows from the expression \eqref{wee}.  The proof is complete.
      \endproof

      For the nonlinear programming,  the mapping $h$ in definition of   $C^2$-conic reducible sets  can be chosen to be an affine mapping and thus $\widetilde{\mathcal{H}}_\lambda u=0$ for all $u\in \R^n.$  In this case, the implication $[{\bf (i)}\Rightarrow {\bf (iv)}]$ was established  by Guo et al. \cite[Theorem 2.1]{GLY13} under the calmness condition, which is weaker than MSCQ.





The following theorem combining Theorem~\ref{thm1ConicP} and the discussion at \eqref{seconic} is the main result of this section, which provides  second order  necessary and sufficient conditions for a stationary point of  \eqref{conicP} to be a strong local minimizer.
\begin{Theorem}\label{thmConicP} {\rm \bf  (second order  characterizations of  strong local minimizer  for  conic program  under MSCQ and $\mathcal{C}^2$-cone reducibility constraint).}\label{main4}  Suppose that  $\bar x$ is a stationary  point of  \eqref{conicP} satisfying MSCQ and that $\Theta$ is $\mathcal{C}^2$-cone reducible at $\oy=h(\ox)$ to a pointed closed convex cone $C$. The following assertions are equivalent:

${\bf (i)}$  $\bar x$ is a strong local minimizer to  problem  \eqref{conicP}.

${\bf (ii)}$ $\bar x$ is a local minimizer to problem  \eqref{conicP} and $\partial f$ is strongly metrically subregular at $\ox$ for $0$.

${\bf (iii)}$ $D(\partial f)(\ox|0)$ is positive definite in the sense of \eqref{SO2}.

${\bf (iv)}$ There exists $\kappa>0$ such that
\begin{equation}\label{SO2C}
\la z,w\ra\geq \kappa\|w\|^2 \quad \mbox{for all}\quad z\in D(\partial f)(\ox|0)(w), w\in \dom D(\partial f)(\ox|0).
\end{equation}

${\bf (v)}$ The sufficient condition of the second-kind in {\rm Definition~\ref{ew}} holds at $\ox$, i.e., there exists $\kappa>0$ such that for each $w\in \dom D(\partial_p f)(\ox|0)$ with $\|w\|=1$ condition \eqref{exi} holds.

${\bf (vi)}$ The {\em second order  sufficient condition} holds at $\ox$ in the sense that for each $w\in \overline{K}\setminus\{0\}$ one has
 $$\max \big\{w^T\nabla^2L(\bar x, \lambda)w+w^T\widetilde{\mathcal{H}}_\lambda w\ |\  \lambda \in  \Lambda(\bar x)\big\}>0.$$

%

 ${\bf (vii)}$ There exists $\kappa>0$ such that
\begin{equation}\label{SOk}
\max \big\{w^T\nabla^2L(\bar x, \lambda)w+ w^T\widetilde{\mathcal{H}}_{\lambda} w\ |\  \lambda \in  \Lambda(\bar x)\big\}\geq \kappa\|w\|^2\quad \mbox{for all}\  w\in \overline{K}.
\end{equation}

 \noindent Moreover,  if one of the assertions ${\bf (i)}-{\bf (vii)}$ holds  then
\begin{eqnarray}\begin{array}{ll}
{\rm QG}((P);\ox)\disp=\inf_{w\in \overline{K},\|w\|=1}\max \big\{w^T\nabla^2L(\bar x, \lambda)w+ w^T\widetilde{\mathcal{H}}_{\lambda} w\ |\  \lambda \in  \Lambda(\bar x)\big\}.\label{QGP}
\end{array}
\end{eqnarray}
  \end{Theorem}
\noindent{\bf Proof.} We see that  $\bar x$ is a strong local minimizer for  \eqref{conicP} with modulus $\kk$  if and only if   $\bar x$ is a local minimizer of the function
$$g_\kk(x):=g(x)-\dfrac{\kk}{2}\|x-\bar x\|^2\quad\mbox{over} \ \Gamma.$$
Applying Theorem~\ref{thm1ConicP} by replacing the function $g$ there by $g_\kk(x)$ along with  using the sum rules of the graphical derivative,  the latter  gives us that
 \begin{equation*}
\la z,w\ra\geq \kappa\|w\|^2 \quad \mbox{for all}\quad z\in D(\partial f)(\ox|0)(w), w\in \dom D(\partial f)(\ox|0),
\end{equation*}
 which verifies  $[{\bf (i)}\Rightarrow {\bf (iv)}]$. By \eqref{wee}, we have $[{\bf (iv)}\Leftrightarrow{\bf (vii)}]$. This allows us to justify the ``$\le$'' inequality in \eqref{QGP} under the validity of ${\bf (i)}$.

Note further from \eqref{norma} and \eqref{sumr} that $\partial_pf(x)=\partial f(x)$ for $x$ around $\ox$. This together with \eqref{wee} shows the equivalence between ${\bf (iv)}$ and ${\bf (v)}$. By \eqref{wee} again, we have $[{\bf (vii)}\Rightarrow {\bf (vi)}\Rightarrow {\bf (iii)}]$. The implication $[{\bf (iii)}\Rightarrow {\bf (ii)}\Rightarrow {\bf (i)}]$ follows from Theorem~\ref{thm1}. This verifies the equivalence of {\bf (i)}--{\bf(vii)}. Finally, the   ``$\ge$'' inequality in \eqref{QGP} is a consequence of \eqref{eq} in Theorem~\ref{thm0}.  The proof is complete. \endproof

The equivalence between {\bf (i)} and {\bf (vi)} above was established under RCQ in Bonnans and Shapiro  \cite[Theorem~1.137]{BS00}. Earlier version of this no gap second order  optimality condition for nonlinear programming was proved in Ioffe \cite{I79} and Ben-Tal \cite{BT80} under the Mangasarian-Fromovitz constraint qualification (MFCQ) recalled later in \eqref{MFCQ}, a particular of RCQ \eqref{RCQ}.  Next let us provide an example, where MSCQ holds while RCQ \eqref{RCQ} does not but all assertions {\bf (i)} - {\bf (vii)} in Theorem~\ref{thmConicP} are satisfied.
\begin{Example}{\rm Consider the problem $(EP)$ as follows
 \begin{equation}\label{E1ICCP} \begin{array}{rl} &\min\limits_{x\in \R^3} \quad g(x)\quad \mbox{s.t.} \quad q(x)\in \mathcal{Q}_3,
 \end{array}
 \end{equation}
 where $\mathcal{Q}_3=\big\{(s_0,s_1,s_2)\in\R^3\ |\  s_0\geq \sqrt{s_1^2+s_2^2} \big\}$ is the second order  cone in $\R^3,$  and the objective function and the constraint mapping are given, respectively, by
 $$g(x):=\frac{1}{2}x_1^2+x_2^2\,\  \mbox{and}\,\      q(x):=\big(2x_2^2, x_2^2-x_3, x_2^2+x_3\big),$$
 for all $ x=(x_1,x_2, x_3)\in \R^3.$ It is well-known that the second order  cone $\mathcal{Q}_3$ is $\mathcal{C}^2$-cone reducible. Moreover, define
 $$\Gamma:=\big\{x |\; q(x)\in \mathcal{Q}_3\big\}=\big\{x=(x_1,x_2, x_3)\in \R^3 |\; x_2^2\geq |x_3|\big\}.$$
With   $\bar  x:=(0,0,0)$, we see that
 $$\begin{array}{rl} g(x)-g(\bar x)&=\frac{1}{2}x_1^2+x_2^2\\
 &=\frac{1}{2}x_1^2+\frac{1}{2}x_2^2+\frac{1}{2}x_2^2\\
 &\geq \frac{1}{2}x_1^2+\frac{1}{2}x_2^2+\frac{1}{2}|x_3|\\
 &\geq \frac{1}{2}(x_1^2+x_2^2+x_3^2)\\
 &=\frac{1}{2}\|x\|^2,
 \end{array}
 $$
for all  $x=(x_1,x_2,x_3)\in \Gamma$  with $|x_3|\leq 1.$
 So   $\bar  x$ is a strong local minimizer to problem  \eqref{E1ICCP}.  

It is easy to check that
 $$\nabla q(\bar x)=\begin{pmatrix} 0&0& 0\\
 0&0&-1\\
 0&0& 1
 \end{pmatrix},\ N_{\mathcal{Q}_3}\big(q(\bar x)\big)=-\mathcal{Q}_3.$$
Hence we have
 $$N_{\mathcal{Q}_3}\big(q(\bar x)\big)\cap {\rm ker}\nabla q(\bar x)^T=\bigcup\limits_{t\in \R}(-\infty, -\sqrt{2}|t|]\times\{t\}\times \{t\},$$
which shows that RCQ \eqref{RCQ} is not satisfied at $\bar x.$

Next let us verify the validity of MSCQ at $\bar x.$ Since $q(x)\in -\mathcal{Q}_3$ if and only if $x_2=x_3=0$, we have $d\big(q(x); \mathcal{Q}_3\big)=0$ when $q(x)\in \mathcal{Q}_3\cup (-\mathcal{Q}_3)$. When $q(x)\notin \mathcal{Q}_3\cup (-\mathcal{Q}_3)$, note that
\[
\sqrt{(x_2^2-x_3)^2+(x_2^2+x_3)^2}-2x_2^2=\sqrt{2x_2^4+2x_3^2}-2x_2^2\geq (|x_3|+x_2^2)-  2x_2^2=|x_3|-x_2^2>0.
\]
It follows that
$$d\big(q(x); \mathcal{Q}_3\big)=\begin{cases} \quad\quad 0\quad\quad\quad \quad \quad\quad \quad \quad\quad  \mbox{if}\  q(x)\in \mathcal{Q}_3\cup (-\mathcal{Q}_3),\\
\frac{1}{\sqrt{2}}\big(\sqrt{2x_2^4+2x_3^2}-2x_2^2\big)\quad\  \mbox{if}\  q(x)\not\in \mathcal{Q}_3\cup (-\mathcal{Q}_3).
\end{cases}$$
When  $q(x)\in \mathcal{Q}_3\cup (-\mathcal{Q}_3)$, we have $d(x;\Gamma)=d\big(q(x); \mathcal{Q}_3\big)=0.$ When $q(x)\not\in \mathcal{Q}_3\cup (-\mathcal{Q}_3)$, define  $u:=\big(x_1, x_2, \frac{x_3}{|x_3|} x_2^2\big)\in \Gamma,$  snd observe that
$$d(x;\Gamma)\leq \|x-u\|=\big|x_3- \frac{x_3}{|x_3|} x_2^2\big|=\big| |x_3|-x_2^2\big|= |x_3|-x_2^2\leq \sqrt{2}d\big(q(x), \mathcal{Q}_3\big).$$
This shows that  MSCQ holds at $\bar x.$ By Theorem~\ref{thmConicP}, all  the assertions {\bf (i)}-{\bf (vii)}   hold.
   }
\end{Example}

Next let us consider a specific form of problem \eqref{conicP}, which is a standard nonlinear programming problem:
  \begin{equation}\label{tiltMP}
  \begin{array}{rl} &\disp\min_{x\in \R^n} \quad g(x)\quad \mbox{subject to} \quad q(x):=(q_1(x), \ldots,q_m(x))\in \Theta:=\R^m_-,
 \end{array}
 \end{equation}
where $g:\R^n\rightarrow\R$ and $q_i:\R^n\rightarrow\R$, $i=1\ldots,m$  are twice continuously differentiable functions. In this case,  by \cite[Example 3.139]{BS00},  $\R^m_-$ is a $\mathcal{C}^2$-conic reducible set and the function $h$ in \eqref{hV}  is chosen as an affine mapping. Thus $\widetilde{\mathcal{H}}_\lambda w=0$ for all $w\in \R^n$ in \eqref{Hti}.

We say  the {\em constant rank constraint qualification} (CRCQ, in brief) holds at $\ox$ if there is a neighborhood $\mathcal{U}$ of $\ox$ such that the gradient system $\{\nabla q_i(x)|\;i\in J\}$ has the same rank in $\mathcal{U}$ for any index $J\subset I(\ox):=\big\{i\in \{1,\ldots,m\}|\; q_i(\ox)=0\big\}$. It is worth noting here that CRCQ  implies  MSCQ; see, e.g.,  \cite{MS11}. Moreover, CRCQ is independent from the the conventional MFCQ at $\ox$, a variant of RCQ \eqref{RCQ} on nonlinear programming:
\begin{equation}\label{MFCQ}
\exists\, d\in \R^n: \la \nabla q_i(\ox),d\ra<0\quad\mbox{for all}\quad  i\in I(\ox).
\end{equation}

According to \cite{PR2,MR12}, a point  $\ox\in\Gamma$ is called   a {\it tilt-stable local minimizer} to problem \eqref{tiltMP}  if there is a number $\gg>0$ such that for $v\in \R^n$ the optimal solution set denoted by $M(v)$ to the following perturbed problem:
\begin{equation}\label{2.10}
  \begin{array}{rl} &\disp\min_{x\in \B_\gg(\ox)} \quad g(x)-\la v,x\ra\quad \mbox{subject to} \quad q(x)\in \R^m_-
 \end{array}
 \end{equation}
is single-valued and Lipschitz continuous  on some neighborhood of $0\in \R^n$ with $M_\gg(0)=\ox$. It is known from \cite[Theorem~3.2]{MN15} that tilt stability could be characterized by a {\em uniform} version of the quadratic growth condition \eqref{GC}. Under CRCQ, we show next that $\ox$ is a tilt stable minimizer to \eqref{tiltMP} if and only if it is just a strong local minimizer, i.e., the pointbased quadratic growth condition \eqref{GC} holds at $\ox$.

\begin{Corollary} \label{corTilt-Str} {\bf (Tilt-stable minimizer and strong local minimizer under CRCQ)} For nonlinear programming problem \eqref{tiltMP}, suppose that the constraint system of  $q(x)\in \R^m_-$   satisfies CRCQ at a stationary point $\bar x$. Then $\bar x$  is a strong local minimizer    if and only if it is a tilt-stable minimizer.
\end{Corollary}
\noindent{\bf Proof.} Suppose that $\bar x$  is a strong local minimizer. Then, by Theorem~\ref{thmConicP}, there exists $\kappa>0$ such that
\begin{equation}\label{mi2}
\max \big\{w^T\nabla^2L(\bar x, \lambda)w\ |\  \lambda \in \Lambda(\bar x)\big\}\geq \kappa\|w\|^2 \quad \mbox{for all}\  w\in \overline{K}.
\end{equation}
  Since CRCQ holds at $\bar x$, by   \cite[Proposition 5.3]{GM15},   $ \overline{\Lambda}(\ox;w)=\Lambda(\bar x)$. For each $\bar\lambda\in \overline{\Lambda}(\ox;w)$, note that
  $$w^T\nabla^2L(\bar x, \bar \lambda)w=\max \big\{w^T\nabla^2L(\bar x, \lambda)w\ |\  \lambda \in \Lambda(\bar x)\big\}.$$
This together with \eqref{mi2} tells us that
    $$w^T\nabla^2L(\bar x, \lambda)w\geq \kappa\|w\|^2 \quad \mbox{for all}\  w\in \overline{K},\ \lambda \in \Lambda(\bar x).$$
 By \cite[Theorem 7.7]{GM15},   $\bar x$ is a tilt-stable local minimizer .

Conversely, if $\bar x$ is a tilt-stable local minimizer, then from \cite[Theorem 3.3]{CHN18} it follows that
\begin{equation*}
\la z,w\ra>0 \quad \mbox{for all}\quad z\in D(\partial f)(\ox|0)(w), w\neq 0.
\end{equation*}
By Theorem~\ref{thmConicP}, $\bar x$ is a strong local minimizer.   The proof is complete. \endproof

The following example taken from \cite[Example 8.5]{GM15} shows that the conclusion of Corollary~\ref{corTilt-Str}  fails if CRCQ is replaced by  MFCQ \eqref{MFCQ}.
\begin{Example}{\rm Consider the following optimization problem
 \begin{equation*}\label{E1NLP} \begin{array}{rl} &\min\limits_{x\in \R^3} -x_1+\frac{1}{2}x_2^2\quad \mbox{subject to} \quad x_1-x_2^4+x_3^2\leq 0,\ x_1\leq 0.
 \end{array}
  \end{equation*}
Then $\bar x=(0,0,0)$ is a strong local minimizer to this problem, but it is not tilt stable. It is easy to check that  MFCQ holds at $\bar x$ and  assertions ${\bf (ii)}-{\bf(vii)}$ in Theorem~\ref{thmConicP} are valid.
}
\end{Example}

To end this section, we provide an example of nonlinear problem satisfying  MSCQ   and  all the assertions {\bf (i)}-{\bf(vii)} in Theorem~\ref{thmConicP} holds at $\bar x$, while  $\bar x$ is not a tilt-stable local minimizer. Both MFCQ and CRCQ are not valid at $\ox$ in this problem.

\begin{Example}{\rm  Let us consider the following nonlinear problem taken from \cite[Example 8.3]{GM15}:
 \begin{equation}\label{E2NLP} \begin{array}{rl} &\min\limits_{x\in \R^3} g(x)=-x_1+\frac{1}{2}x_2^2+\frac{1}{2}x_3^2\\
 &\quad \mbox{s.t.} \quad q_1(x)=  x_1-\frac{1}{2}x_2^2\leq 0,\\
 &\quad\quad\quad \  q_2(x)=x_1-\frac{1}{2}x_3^2\leq 0,\\
 &\quad\quad\quad \  q_3(x)=-x_1-\frac{1}{2}x_2^2-\frac{1}{2}x_3^2\leq 0.
 \end{array}
  \end{equation}
  By \cite[Example 8.3]{GM15}, we see that both MFCQ and CRCQ do not hold at $\bar x:=(0,0,0),$ and $\bar x$ is not a tilt-stable local minimizer.  Put
  $$\Gamma=\{x=(x_2,x_2,x_3)\in \R^3 | q_i(x)\leq 0,\ i=1,2,3\}.$$
  We now show that $\bar x$ is a strong local minimizer. Observe that
  \begin{equation}\label{Ex4.7-eq1} g(x)-\frac{1}{8}\|x\|^2=-x_1-\frac{1}{8}x_1^2+\frac{3}{8}x_2^2+\frac{3}{8}x_3^2,\ x=(x_1,x_2,x_3)\in \R^3.\end{equation}
    If $-\frac{1}{8}\leq x_1\leq 0$ then \eqref{Ex4.7-eq1} implies that
    $$g(x)-\frac{1}{8}\|x\|^2\geq 0.$$
   If $x_1 \geq 0,$ $-2\leq x_2\leq 2$  and $x\in \Gamma$  then
$$\begin{cases}\frac{1}{2}x_2^2\geq x_1\geq 0,\\
\frac{1}{2}x_3^2\geq x_1\geq 0,
\end{cases}
$$
which  infers  that
$$\begin{cases} -\frac{1}{2}x_1\geq -\frac{1}{4}x_2^2,\\
-\frac{1}{2}x_1\geq -\frac{1}{4}x_3^2,\\
-\frac{1}{8}x_1^2\geq -\frac{1}{32}x_2^4.
\end{cases}
$$
So, by \eqref{Ex4.7-eq1} we have
$$ \begin{array}{rl}g(x)-\frac{1}{8}\|x\|^2&\geq -\frac{1}{4}x_2^2 -\frac{1}{4}x_3^2-\frac{1}{32}x_2^4+\frac{3}{8}x_2^2+\frac{3}{8}x_3^2\\ \cr
&=\frac{1}{8}x_2^2(1-\frac{1}{4}x_2^2)+\frac{1}{8}x_3^2\geq 0,\end{array}$$
provided that $x\in \Gamma$ with $x_1 \geq 0$ and $-2\leq x_2\leq 2.$
Hence,
$$g(x)-\frac{1}{8}\|x\|^2\geq 0\ \mbox{for all}\  x\in \Gamma\cap\big([-\frac{1}{8},+\infty)\times [-2,2]\times \R\big).$$
Thus $\bar x$ is a strong local minimizer. By Theorem~\ref{thmConicP}, all  the assertions {\bf (i)}-{\bf (vii)}   hold.
}
\end{Example}


\section{Conclusion}
We have used the subgradient graphical derivative to study the quadratic growth condition and strong local minimizer for optimization problems. For the nonsmooth unconstrained problems, we show that the positive definiteness of the subgradient graphical derivative at the proximal stationary point is sufficient for the quadratic growth condition. It becomes necessary in the case of either subdifferentially continuous, prox-regular, twice epi-differentiable or  variational convex functions. Our approach is mainly based on the property of strong metric subregularity on the subdifferential.  For the smooth $\mathcal{C}^2$-reducible cone constrained problems, we obtain full characterizations for strong local minimizers under the metric subregularity constraint qualification, which is strictly weaker than the classical Robinson's constraint qualification. One of our main results shows that the traditional second order sufficient  condition could characterize the quadratic growth condition under MSCQ. In the future, we intend  to extend our applications to different optimization problems including composite functions \cite{BCS99,BS00,RW98} and mathematical programs with equilibrium constraints \cite{GY17,GLY13}. This extension will require further computation on the second order  structures of subgradient graphical derivative presented in our paper. Another direction that cause our intention is to study the epi-differentiability of composite functions and conic constrained optimization problems under metric subregularity  conditions in order to use the advantage of Theorem~\ref{thm3}.

\end{document}